\documentclass[final,leqno]{siamltex704}

\usepackage{amsmath,amssymb,amsfonts}
\usepackage{color,url,graphicx,epsfig,subfigure}
\usepackage{cancel,enumerate,version}
\newtheorem{remark}[theorem]{Remark}

\newtheorem{example}[theorem]{Example}
\newtheorem{conjecture}[theorem]{Conjecture}
\newcommand{\ml}{\left[ \begin{array}{cccccccccccccccccccccccccc}}
\newcommand{\mr}{\end{array} \right]}
\newcommand{\mll}{\left[ \begin{array}{c|c|c|c|c}}
\newcommand{\yf}{\operatorname{fvct}}

\newcommand{\until}[1]{\{1,\dots, #1\}}
\newcommand{\subscr}[2]{#1_{\textup{#2}}}

\newcommand\oprocendsymbol{\hbox{$\square$}}
\newcommand\oprocend{\relax\ifmmode\else\unskip\hfill\fi\oprocendsymbol}

\renewcommand{\theenumi}{(\roman{enumi})}

\newcommand{\real}{\mathbb{R}}
\begin{document}

\title{Opinion Dynamics in Heterogeneous Networks: Convergence Conjectures
  and Theorems\thanks{This work was supported in part by the UCSB Institute
    for Collaborative Biotechnology through grant DAAD19-03-D004 from the
    U.S. Army Research Office. This article was submitted to the SIAM
    Journal on Control and Optimization in March 2011.}}

\author{Anahita~Mirtabatabaei \hspace{0.5in} Francesco~Bullo\thanks{Anahita
    Mirtabatabaei and Francesco Bullo are with the Center for Control,
    Dynamical Systems, and Computation, University of California, Santa
    Barbara, Santa Barbara, CA 93106, USA,\tt{
      \{mirtabatabaei,bullo\}@engineering.ucsb.edu}}} \date{\today}
\maketitle

\begin{abstract}
  Recently, significant attention has been dedicated to the models of opinion
  dynamics in which opinions are described by real numbers, and agents
  update their opinions synchronously by averaging their neighbors'
  opinions.  The neighbors of each agent can be defined as either (1) those
  agents whose opinions are in its ``confidence range,'' or (2) those
  agents whose ``influence range'' contain the agent's opinion.  The former
  definition is employed in Hegselmann and Krause's bounded confidence
  model, and the latter  is novel here. As the
  confidence and influence ranges are distinct for each agent, the
  heterogeneous state-dependent interconnection topology leads to a
  poorly-understood complex dynamic behavior.  In both models, we classify
  the agents via their interconnection topology and, accordingly,
  compute the equilibria of the system.
  Then, we define a positive invariant set centered
  at each equilibrium opinion vector.  We show that if a trajectory enters
  one such set, then it converges to a steady state with constant
  interconnection topology. This result gives us a novel sufficient
  condition for both models to establish convergence, and is consistent
  with our conjecture that all trajectories of the bounded confidence and
  influence models eventually converge to a steady state under fixed
  topology.
\end{abstract} 


\begin{AMS} 
 91D30, % Social networks
 93A14, % Decentralized systems
 37A99,  % Dynamical systems and ergodic theory: None of the above, but in this section
 91C20, % Clustering
 37B55, % Nonautonomous dynamical systems
 91B14 %Social choice
\end{AMS}

\section{Introduction}
Models of social networks are structures made up of individuals that are
tied based on their interdependency. Such models explain the confidence or
influence flow in populations without relying on detailed social
psychological findings.  The process of opinion dynamics evolves along the
networks of social confidence or influence and affects the structure of the
network itself.
A common feature among
many models of opinion dynamics is bounded confidence or influence, which
means that an individual only interacts with those whose opinions are close
enough to its own. This idea reflects the psychological concept called
\emph{selective exposure} \cite{WM-FC-RS:08}. Broadly defined, ``selective
exposure refers to behaviors that bring the communication content within
reach of one's sensory apparatus'' \cite{DZ-JB:85}.
In the field of social networks, opinion dynamics is of high interest in
many areas including: politics, as in voting prediction \cite{EBN:05};
physics, as in spinning particles \cite{EN-PK-FV-SR:03}; sociology, as in
the diffusion of innovation \cite{TV:96}, the electronic exchange of personal
information \cite{EM-NL-KF-RP-HL-MM-SE-XZ:08}, and language change
\cite{JT-MT-AN:07,FC-SS:07}; and finally economics, as in price
change \cite{KS-RW:02}.

\subsection{Literature Review} 
The study of opinion dynamics and social networks goes back to the early
work by J.R.P.\ French \cite{JF:56} on ``A Formal Theory of Social Power.''
This work explores the patterns of interpersonal relations and agreements
that can explain the influence process in groups of agents. Subsequently,
F.\ Harary provides a necessary and sufficient condition to reach a
consensus in French's model of power networks \cite{FH:59}. The modeling of
``continuous opinion dynamics'', in which opinions are represented by real
positive numbers, is initially studied in \cite{MS:61,SC-ES:77,CJ-EJ:86}.
In contrast to the classical case of ``binary opinion dynamics''
\cite{GW-GB:02,SG:97}, the continuous case deals with the problem of what
happens to the worthiness of a choice or the probability of choosing one
decision over another.  Recently, \emph{bounded confidence} (BC) models of
opinion dynamics, a label coined by Krause in 1998 \cite{UK:00}, have
received significant attention. BC models are models of continuous opinion
dynamics in which agents have bounded confidence in others opinions. The
first version of BC models was formulated by Hegselmann and Krause
\cite{RH-UK:02}, called HK model, where agents synchronously update their
opinions by averaging all opinions in their confidence bound. The other
popular version of BC models was developed and investigated by Deffuant and
Weisbuch \cite{GW-GD-FA-JN:02}, called DW model. The HK and DW models are
very similar, they differ in their update rule: in DW model a
pairwise-sequential updating procedure is employed instead of the
synchronized one.

In the HK model, the set of neighbors of the $i$th agent is defined as
those agents whose opinions differ from the $i$th opinion by less than the
$i$th confidence bound. Hence, this model is dealing with
\emph{endogenously} changing topologies, that is, state dependent or
changing from inside, in contrast to the \emph{exogenously} changing
topologies. For instance,
\cite{AJ-JL-ASM:02,LM:05,FC-SS:07,MC-ASM-BDOA:06} study a synchronized
linear averaging model with time-dependent exogenously changing
topologies. The HK models are classified based on various factors: a model
is called agent- or density-based if its number of agents is finite or
infinite, respectively; and a model is called homogeneous or heterogeneous
if its confidence bounds are uniform or agent-dependent, respectively.  The
convergence of both agent- and density-based homogeneous HK models are
discussed in \cite{VDB-JMH-JNT:09}.  The agent-based homogeneous HK system
is proved to reach a fixed state in finite time \cite{JD:01}, the time
complexity of this convergence is discussed in \cite{SM-FB-JC-EF:05m}, its
stabilization theorem is given in \cite{JL:05a}, and its rate of
convergence to a global consensus is studied in \cite{AO-JNT:07}.  The
heterogeneous HK model is studied by Lorenz who reformulated the HK
dynamics as an interactive Markov chain \cite{JL:06b} and  analyzed the
effects of heterogeneous confidence bounds \cite{JL:09}. %
The convergence of the agent-based heterogeneous HK systems is
experimentally observed, but its proof is still an open problem.

\subsection{Contributions} 
In this paper, to distinguish between the HK and DW models, we call a
discrete-time agent-based heterogeneous HK model a \emph{synchronized
  bounded confidence} (\emph{SBC}) model. Additionally, we introduce a
model similar to the SBC model and call it the \emph{synchronized bounded
  influence} (\emph{SBI}) model. The difference is that in an SBI model the
set of neighbors of the $i$th agent is defined as those agents $j$ whose
influence range contain the $i$th agent's opinion. We analyze SBC and SBI
models with heterogeneous bounds of confidence or influence, respectively.
Indeed, if the SBC and SBI models have agents with homogeneous bounds, then
both models are equivalent to the homogeneous HK model.  These
heterogeneous models of opinion dynamics, in spite of the considerable
complexity of their dynamics, describe features of a real society that
cannot be explained by homogeneous models.  Specifically, the behavior of a
heterogeneous opinion dynamics model is richer than that of a homogeneous
model in the following ways: (1) an agent may trust an individual but not
be trusted back by that same individual; (2) in steady state, an agent can
keep its own opinion constant while listening to dissimilar opinions,
(whereas, in the steady state of a homogeneous society, any two agents are
either disconnected or in consensus); (3) one can observe
\emph{pseudo-stable} configurations, i.e., ``configurations that have a
subset of agents that is stationary, and the rest are the reason for
further dynamics'' \cite{DU:03}; (4) it is possible for two disconnected
agents to reconnect, and an agent with a large bound of confidence or
influence can pull clusters of agents towards or away from each other; (5)
the order of opinions is not preserved along the system evolution; (6)
"given an average level of confidence, diversity of bounds of confidence
enhances the chances for consensus" \cite{JL:09}; and (7) convergence in
infinite time is possible.

Based on numerical evidence, we formulate our main conjecture that along
the evolution of an SBC or SBI system there exists a finite time, after
which the topology of the interconnection network remains unchanged, and as
a result, the trajectory converges to a limiting opinion vector.  We also
observe that each trajectory either reaches a fixed state in finite time or
exhibits a pseudo-stable behavior. This observation is verified assuming
that the main conjecture is true.  Furthermore, the following results put
together partly prove our main conjecture: (1) We design an appropriate
classification of agents in both SBC and SBI systems. This classification
is a function of state-dependent interconnection topology of the system,
and can explain the observed pseudo-stable behavior. (2) We introduce the
new notion of \emph{final value at constant topology}, and based on our
classification, we formulate the map under which this value is an image of
the current opinion vector. The set of final values at constant topology is
a superset of the equilibria of the system. We derive necessary and
sufficient conditions for the final value at constant topology to be an
equilibrium vector.  (3) For each equilibrium opinion vector, we define its
\emph{equi-topology neighborhood} and \emph{invariant equi-topology
  neighborhood}. We show that if a trajectory enters the invariant
equi-topology neighborhood of an equilibrium vector, then it remains
confined to its equi-topology neighborhood, and sustains an interconnection
topology equal to that of the equilibrium vector. This fact establishes a
novel and simple sufficient condition under which: the initial opinion
vector converges to a steady state; the topology of the interconnection
network remains unchanged; and the limiting opinion vector is equal to the
final value at constant topology of the initial opinion vector.  (4) We
explore some interesting behavior of classes of agents when they update
their opinions under fixed interconnection topology for infinite time. For
instance, we compute agents rates and directions of convergence, and show
the existence of a leader group for each group of agents that determines
the follower's rate and direction of convergence.

In our extensive simulation results, we observe that for uniformly randomly
generated initial opinion vector and bounds vector, the SBC and SBI
trajectories eventually satisfy our novel sufficient condition for
convergence with probability one. We give some intuitive explanation for
this observation. Finally, we conjecture that the SBI trajectories reach a
fixed state in finite time more often than the SBC trajectories. To
substantiate this conjecture, we present a sufficient condition for SBC and
SBI systems separately that guarantees reaching an \emph{agreement opinion
  vector}, which occurs in finite time and is a general case of reaching a
global consensus. Based on these sufficient conditions, we explain our last
conjecture.

\subsection{Organization}

This paper is organized as follows. In Section~\ref{model}, the
mathematical models, conjectures, agents classification, and spectral
properties of the adjacency matrices are presented. In
Section~\ref{final_value}, the final value at constant topology is
introduced and characterized. Section~\ref{sufficientcond} contains the
novel sufficient condition for constant topology and convergence, and for
convergence in finite time.  In Section~\ref{expsec}, the simulation
results and intuitive explanations are presented. In
Section~\ref{fixedtopbehavior}, the behavior of the system assuming that
its interconnection topology remains unchanged in a long run is
analyzed. Finally, Section~\ref{conclusion} contains the conclusion and
open questions.

\section{Mathematical Models}\label{model}
Consider $n$ interacting agents and assume that each agent's opinion is
expressed by a real number, say $y_i$ for agent $i\in\until{n}$.  In
\emph{bounded confidence interaction}, the opinion $y_i$ is affected by the
opinion $y_j$ if $|y_i-y_j| \le r_i$, where the positive number $r_i$ is the
\emph{confidence bound} of agent $i$.  In \emph{bounded influence
  interaction}, the opinion $y_i$ is affected by the opinion $y_j$ if
$|y_i-y_j| \le r_j$, where the positive number $r_j$ is the \emph{influence
  bound} of agent $j$.  The \emph{opinion vector} $y\in \real^n$ and the
\emph{bounds vector} $r\in\real^n_{>0}$ are obtained by stacking all
$y_i$'s and $r_i$'s, respectively.  We associate to each opinion vector $y$
two digraphs, both with nodes $\until{n}$ and edge set defined as follows:
denoting the set of out-neighbors of node $i$ by $\mathcal{N}_i (y)$
\begin{itemize}
\item in a \emph{synchronized bounded confidence} (SBC) digraph,
  $\mathcal{N}_i (y) = \{ j\in\until{n} : |y_i - y_j| \le r_i \}$; and
\item in a \emph{synchronized bounded influence} (SBI) digraph,
  $\mathcal{N}_i (y) = \{ j\in\until{n} : |y_i - y_j| \le r_j \}$.
\end{itemize}
We let $G_r(y)$ denote one of the two \emph{proximity digraphs}, its
precise meaning being clear from the context.

\begin{figure}
  \centering
  \includegraphics[width=.47\textwidth,keepaspectratio]{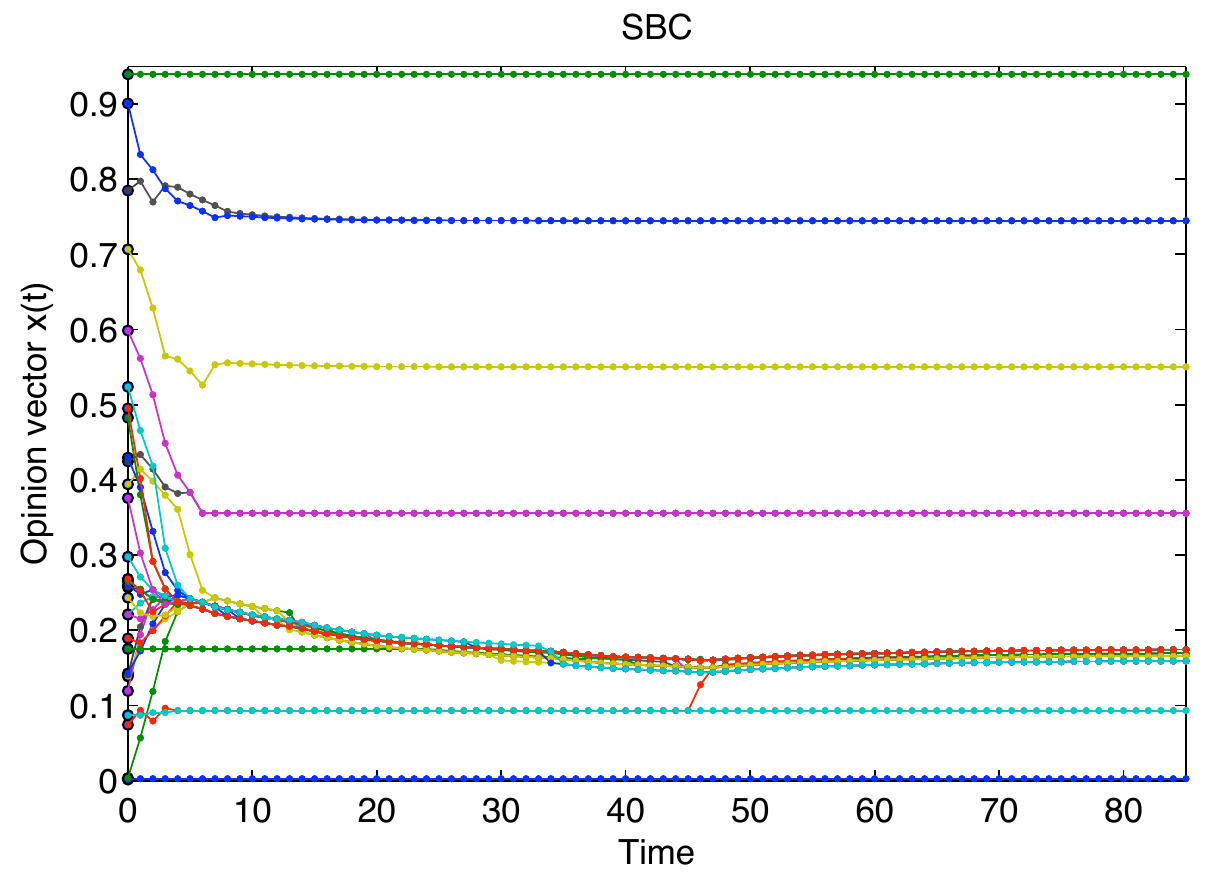}
  \includegraphics[width=.47\textwidth,keepaspectratio]{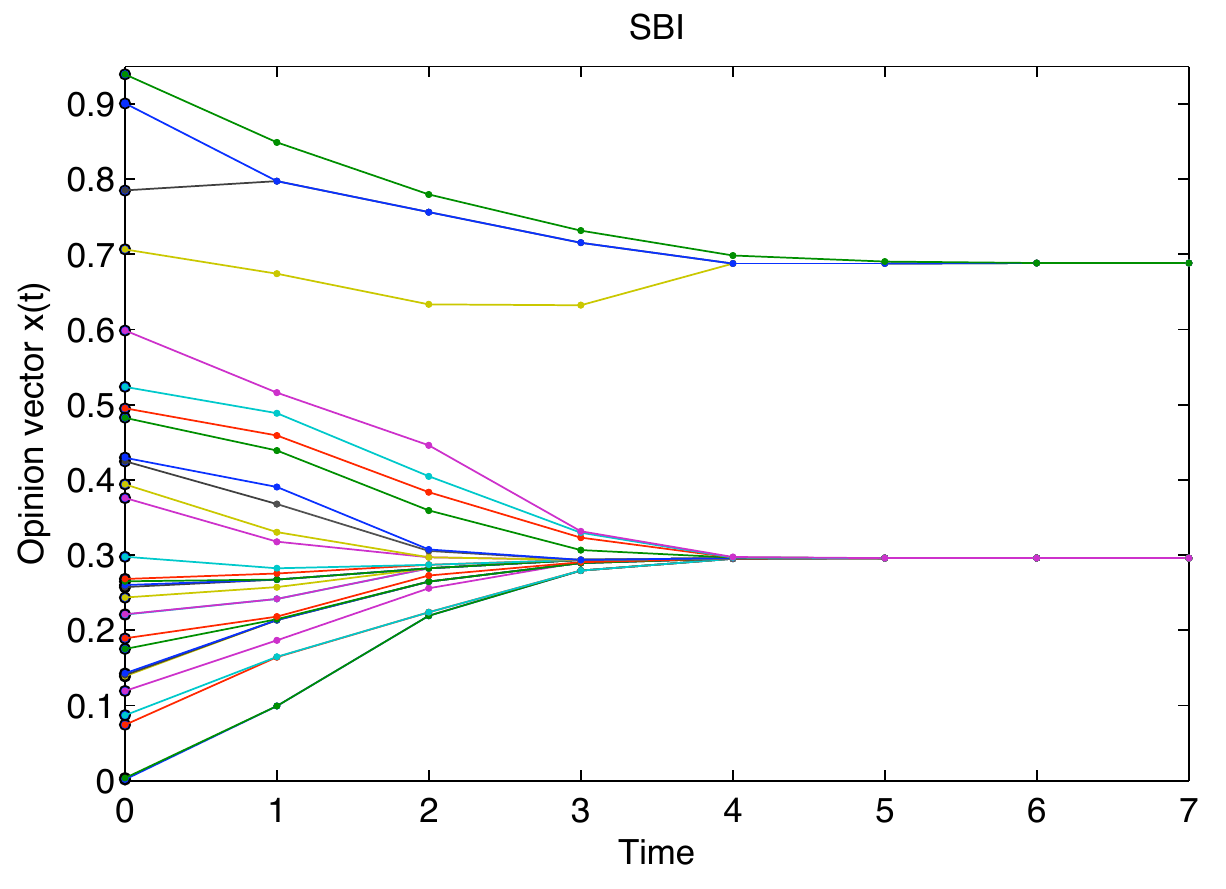}
  \caption{The trajectory of an SBC system (left) and an SBI system (right) are illustrated.
Both systems have the same initial opinion vector and bounds vectors that are randomly generated. However, the SBI
  trajectory reaches a fixed state in six time steps, while the SBC trajectory converges in infinite time.
The interconnection topology of the agents in the SBC system 
 remains constant after $t=64$, and hence its trajectory exhibits a pseudo-stable behavior.}\label{evolutionfig}
\end{figure}

We associate to the SBC and SBI digraphs two dynamical systems, called the
\emph{SBC} and \emph{SBI systems} respectively.  Both dynamical systems
update a trajectory $x: \mathbb{N} \rightarrow \real^n$ according to the
discrete-time and continuous-state rule
\begin{equation} \label{system}
  x(t + 1)  =   A(x(t)) x(t),
\end{equation}
where the $i,j$ entry of the \emph{adjacency matrix} $A(y) \in
\real^{n\times n}$ for any $y \in \real^n$ is defined by
\begin{align*}
  a_{ij} (y) & = 
  \begin{cases} 
   \frac{1}{\displaystyle |\mathcal{N}_i (y)|},  & \mbox{if } j \in \mathcal{N}_i(y),\\ 
    0, & \mbox{if } j \notin \mathcal{N}_i(y),
  \end{cases}
\end{align*}
and $|\mathcal{N}_i (y)|$ is the cardinality of $\mathcal{N}_i (y)$. Note
that $i \in \mathcal{N}_i (y)$, in other words, every agent has some
self-confidence or self-influence. This assumption is a key factor in the
convergence of infinite products of adjacency matrices \cite{JL:06}.  In the
following, we present our conjectures on SBC and SBI
systems, and the trajectories of Figure~\ref{evolutionfig} support these
conjectures.
\begin{conjecture}[Existence of a limiting opinion vector]\label{xinfconj}
  Every trajectory of an SBC or SBI system converges to a limiting opinion
  vector.
\end{conjecture}
\begin{conjecture}[Constant-topology in finite time]\label{conjtau}
  For any trajectory $x(t)$ of an SBC or SBI system, there exists a finite
  time $\tau$ after which the state-dependent interconnection topology, or
  equivalently $G_r(x(t))$, remains constant.
\end{conjecture}

Before proceeding, let us define a term borrowed from \cite{DU:03}. A trajectory $x(t)\in\real^n$ that is converging to limiting opinion vector $x_\infty\in\real^n$ is said to have a \emph{pseudo-stable} behavior after $\tau$, if the node set $\mathcal{V} = \until{n}$ is composed of two non-empty subsets $\mathcal{V}_{\text{fixed}}$ and $\mathcal{V}_{\text{converging}}$ such that, for all $t\ge \tau$,
\begin{equation}\label{pseudoeq} 
    \begin{cases} 
      x_i(t) = x_{\infty,i}, & \mbox{if } \; i \in \mathcal{V}_{\text{fixed}},\\
      x_i(t) < x_i(t+1) < x_{\infty,i} \;\; or \;\; x_i(t) > x_i(t+1) >
      x_{\infty,i} , & \mbox{if } \; i \in \mathcal{V}_{\text{converging}}.
    \end{cases}
\end{equation}
\begin{conjecture}[Pseudo-stable behavior]\label{pseudoconj}
  For any SBC or SBI trajectory, there exists a finite time after which the
  trajectory either reaches a fixed state or exhibits a pseudo-stable
  behavior.
\end{conjecture}
\begin{conjecture}[Convergence of SBI systems versus SBC systems]\label{inffiniteconj}
For any initial opinion vector and 
bounds vector that are generated uniformly randomly, the SBI system is more likely to converge in finite time than the SBC system.
\end{conjecture}

This paper aims to study these conjectures. 

\subsection{Agents Classification}\label{class}

In this section, we introduce a classification of agents for both SBC and SBI systems based on their state-dependent interaction topology at each time step. This classification is used later to find the limiting opinion vector and explain the
pseudo-stable behavior.  
First, let us quote some relevant definitions from graph theory,
e.g. see \cite{FB-JC-SM:09}.  A node of a digraph is \emph{globally
  reachable} if it can be reached from any other node by traversing a
directed path. A digraph is \emph{strongly connected} if every node is
globally reachable. A digraph is \emph{weakly connected} if replacing all
of its directed edges with undirected edges produces a connected undirected
graph.  A maximal subgraph which is strongly or weakly connected forms a
\emph{strongly connected component} (SCC) or a \emph{weakly connected
  component} (WCC), respectively. Every digraph $G$ can be decomposed into
either its SCC's or WCC's. Accordingly, the \emph{condensation digraph} of
$G$, denoted $C (G)$, is defined as follows: the nodes of $C (G)$ are the
SCC's of $G$, and there exists a directed edge in $C (G)$ from node $H_1$
to node $H_2$ if and only if there exists a directed edge in $G$ from a
node of $H_1$ to a node of $H_2$. A node with out-degree zero is named a
\emph{sink}. Knowing that the condensation digraphs are acyclic, each WCC
of $C(G)$ is also acyclic and thus has at least one sink. In a digraph, $i$ is a
\emph{predecessor} of $j$ and $j$ is a \emph{successor} of $i$ if 
there exists a directed path from node $i$ to node $j$.

\begin{figure}
  \centering
  \includegraphics[width=3.4in,keepaspectratio]{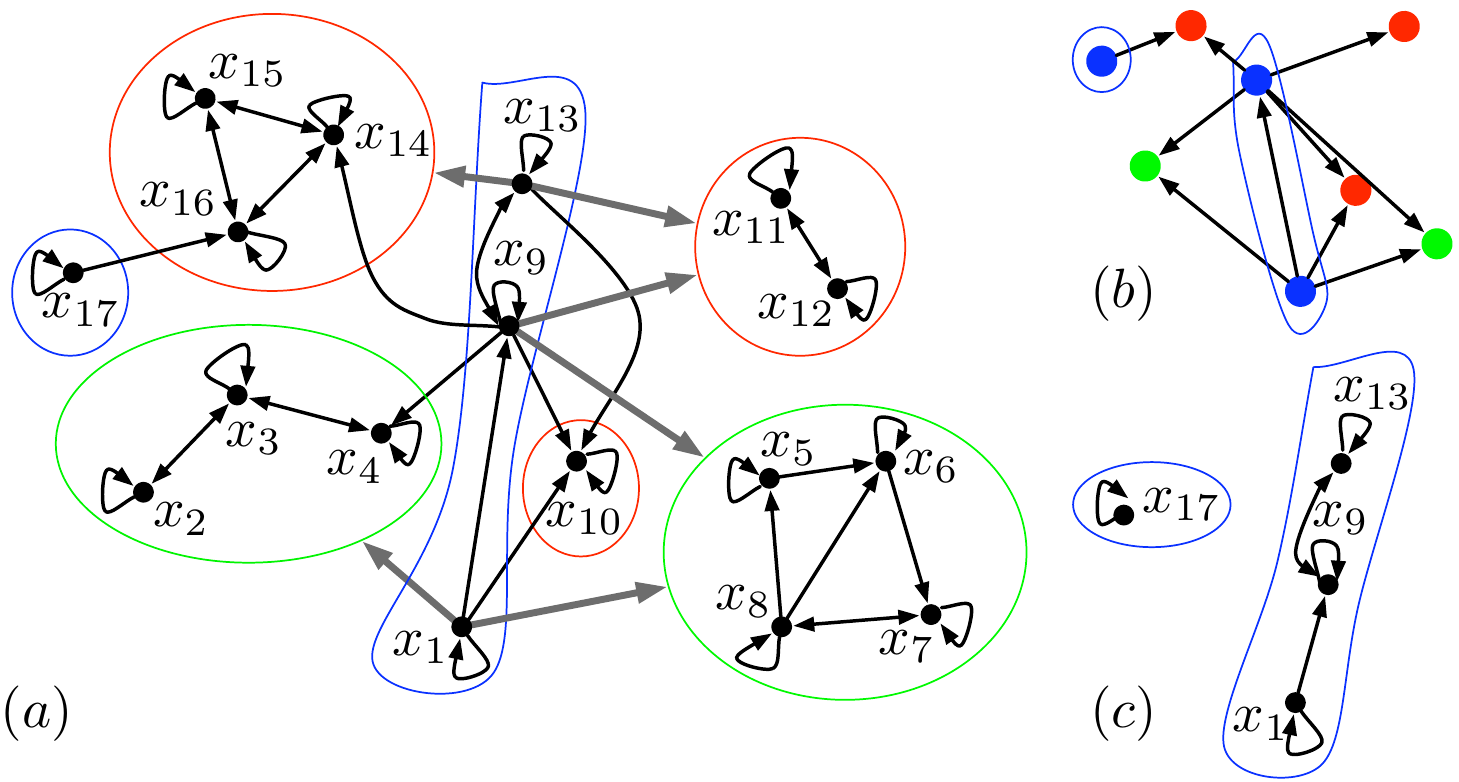}
  \caption{Consider the opinion vector $x=$ [0.1 0.24 0.27 0.3 0.34 0.37
    0.39 0.4 0.5 0.6 0.67 0.68 0.75 0.85 0.86 0.87 1]$^T$ and bounds vector
    $r=$[0.5 0.04 0.04 0.04 0.031 0.021 0.011 0.061 0.25 0.01 0.04 0.03 0.3
    0.07 0.07 0.07 0.135]$^T$: (a) shows the SBC digraph of $x$, $G_r(x)$,
    with its closed- (red), moderate- (green), and open-minded (blue)
    components, and each thick gray edge represents multiple edges to all
    agents in one component; (b) shows the condensation digraph of
    $G_r(x)$; and (c) shows the open-minded subgraph of $G_r(x)$ that is
    composed of two open-minded WCC's.}\label{proximityfig}
\end{figure}

For opinion vector $y \in \real^n$, let $G_r(y)$ denote either of its SBC or SBI digraphs. We classify the SCC's of $G_r(y)$ into three classes.
An SCC of $G_r(y)$ is called a \emph{closed-minded component} if it is a complete subgraph of $G_r(y)$ and corresponds to a sink of $C(G_r(y))$. 
An SCC of $G_r(y)$ is called a \emph{moderate-minded component} if it is a non-complete subgraph of $G_r(y)$ and corresponds to a sink of $C(G_r(y))$. 
The rest of SCC's of $G_r(y)$ are called \emph{open-minded SCC's}. 
Now, the \emph{open-minded subgraph} of $G_r(y)$ is the remaining subgraph after removing $G_r(y)$'s closed- and moderate-minded components and their edges. A WCC of the open-minded subgraph of $G_r(y)$ will be called an \emph{open-minded WCC}, see Figure~\ref{proximityfig}. 
\begin{remark}
Previously, (Lorenz, 2006) classified the agents of an SBC system into two classes of \emph{essential} and \emph{inessential}. An agent is
essential if any of its successors is also a predecessor, and an agent is inessential if it has a successor who is not a predecessor \cite{JL:06}. 
This classification is similar to the one used for Markov chains \cite[Chapter 1.2]{ES:81}. It is easy to see that the closed- and moderate-minded components are in essential class, and the open-minded components are inessential. 
\end{remark}

\subsection{Spectral Properties of Adjacency Matrix}\label{blocks}
For any opinion vector $y \in \real^n$ in an SBC or SBI system \eqref{system}, the adjacency matrix $A(y)$ is a non-negative row-stochastic matrix, and its nonzero diagonal establishes its aperiodicity. 
Since $C(G_r(y))$ is an acyclic digraph, its adjacency matrix is lower-triangular in an appropriate ordering \cite{FB-JC-SM:09}. In such ordering, the adjacency matrix of $G_r(y)$ is lower block triangular.
Based on the classification of the SCC's in $G_r(y)$, we put $A(y)$ into the \emph{canonical form} $\overline{A}(y)$, by an appropriate
 \emph{canonical permutation matrix} $P(y)$, 
\begin{equation}\label{canonical}
\overline{A}(y) = P(y) A(y) P^T(y) = \ml C(y) & 0 & 0\\0 & M(y) & 0\\  \Theta_C(y) & \Theta_M(y)&  \Theta(y) \mr.
\end{equation}
The submatrices $C(y)$, $M(y)$, and $\Theta(y)$ are block diagonal. Each diagonal block $C_i(y)$, with size $n_i(y)$, is the adjacency matrix of the $i$th closed-minded component, and is equal to $C_i(y)=  {\bf 1}_{n_i(y)}{\bf 1}^T_{n_i(y)} / n_i(y).$ Let us call a matrix with such structure a \emph{complete consensus matrix}, whose spectrum is found to be $\{ 1, 0, \dots,0 \}.$ 
Similarly, each diagonal block $M_i(y)$ is the adjacency matrix of the $i$th moderate-minded component. 
Each entry in $\Theta_C(y) $ or $ \Theta_M(y)$ represents an edge from an open-minded node to a closed- or moderate-minded node, respectively. 
Finally, in the submatrix $\Theta(y)$, each diagonal block $\Theta_i(y)$ corresponds to one open-minded WCC, and is block lower triangular and strictly row-substochastic. By strictly row-substochastic we mean a square matrix with nonnegative entries so that every row adds up to at most one, and there exists at least one row whose sum is strictly less than one. 
Note that the adjacency matrix of each SCC in $G_r(y)$ is a diagonal block of $\overline{A}(y)$ and is row-stochastic, nonnegative and primitive. On account of the properties of the open-minded class, the following lemma is proved.
\begin{lemma}\label{rowsum}
For any row $k$ of the submatrix $\Theta(y)$, there exists $p_k \in \mathbb{N}$ such that the $k$th row sum of $\Theta(y)^{p_k}$ is strictly less than 1.
\end{lemma}
\begin{proof}
Every WCC of $C(G_r(y))$ contains at least one sink.
Hence, from any open-minded agent $k$, there exists a directed path of length $p_k$ to an agent $s$ in either a closed- or moderate-minded component. Now, consider the canonical adjacency matrix to the power~$p_k$,
$$\overline{A}(y)^{p_k} = \ml C(y)^{p_k} & 0 & 0\\0 & M(y)^{p_k} & 0\\ \Theta_C^{(p_k)}(y)  & \Theta_M^{(p_k)}(y) &  \Theta(y)^{p_k}\mr.$$ 
Existence of such directed path, by \cite[Lemma 1.32]{FB-JC-SM:09}, implies that the $(k,s)$ entry of  $\overline{A}(y)^{p_k}$, which belongs to either of the submatrices $\Theta_C^{(p_k)}(y)$ or $\Theta_M^{(p_k)}(y)$, is nonzero. Consequently, the $k$th row sum of $\Theta(y)^{p_k}$ is strictly less than 1.
\end{proof}

It follows from Lemma~\ref{rowsum} that $\lim_{t \rightarrow \infty} \Theta(y)^t = {\bf 0}, $
for a proof of which refer to \cite[Theorem 4.3]{ES:81}. 
Therefore, the spectral radius of $\Theta(y)$ is strictly less than one.

\begin{example}\label{exampleblocks}    
Consider the SBC system of Figure~\ref{proximityfig} with the permuted opinion vector $x = [ x_2 \;\;  x_3 \;\;  x_4 \;\;  x_5 \;\;  x_6 \;\;  x_7 \;\;  x_8\;\; x_{10} \;\;  x_{11} \;\;  x_{12} \;\;  x_{14} \;\;  x_{15} \;\;  x_{16} \;\;  x_{17} \;\;  x_9 \;\;  x_1 ]^T$.
Then, the canonical form of the adjacency matrix $A(x)$ contains the following submatrices:

$$C(x) =  \ml 
         1     &&&&&   \\
       &       \frac{1}{2} &   \frac{1}{2}   &&&\\
       &       \frac{1}{2} &   \frac{1}{2}   &&&\\       
         &&&      \frac{1}{3} &   \frac{1}{3}   &   \frac{1}{3} \\    
         &&&      \frac{1}{3} &   \frac{1}{3}   &   \frac{1}{3} \\    
         &&&      \frac{1}{3} &   \frac{1}{3}   &   \frac{1}{3} \mr, \;\;
M(x) =  \ml 
       \frac{1}{2} &   \frac{1}{2}   &&&&&\\
       \frac{1}{3} &   \frac{1}{3}   &   \frac{1}{3} &&&&\\
       &      \frac{1}{2} &   \frac{1}{2}   && &&\\
       &&&       \frac{1}{2} &   \frac{1}{2}   &&\\
       &&&&       \frac{1}{2} &   \frac{1}{2}   &\\
       &&&&&      \frac{1}{2} &   \frac{1}{2}   \\
       &&&       \frac{1}{4} &   \frac{1}{4}   &   \frac{1}{4} &   \frac{1}{4} \mr,$$

$$\Theta(x) = \ml 
  \frac{1}{2}   &&  &\\
  & \frac{1}{8} & \frac{1}{8}&\\
  & \frac{1}{11} &   \frac{1}{11}   &\\
 && \frac{1}{10} &   \frac{1}{10}  \mr, \;\;
\Theta_C(x)  =  \ml
  &&&&& \frac{1}{2} \\
 \frac{1}{8} & \frac{1}{8} &  \frac{1}{8} &   \frac{1}{8} &  \frac{1}{8} &   \frac{1}{8} \\ 
 \frac{1}{11} &  \frac{1}{11} &   \frac{1}{11} &\frac{1}{11} &&\\
 \frac{1}{10} &&&&& \mr, $$

$$\Theta_M(x)  =  \ml
\tiny 0 &&& \dots &&&  0\\
 0 &&& \dots &&&  0 \\
 &   &  &   \frac{1}{11} &  \frac{1}{11} &   \frac{1}{11} &  \frac{1}{11} \\
 \frac{1}{10} &   \frac{1}{10} &  \frac{1}{10} &   \frac{1}{10} &  \frac{1}{10} &   \frac{1}{10} &  \frac{1}{10} \mr.$$
\end{example}

\section{Equilibria and Final Value at Constant Topology}\label{final_value}
An opinion vector $y_0$ is an \emph{equilibrium opinion vector} of the dynamical system~\eqref{system}
if and only if $y_0$ is an eigenvector of the adjacency matrix $A(y_0)$ for the eigenvalue one or, equivalently, $y_0 = A(y_0) y_0$. Next, based on Conjecture~\ref{conjtau}, we introduce the following definition.
\begin{definition}[Final value at constant topology]
For any opinion vector $y\in\real^n$ we define its \emph{final value at constant topology} $\yf : \real^n \rightarrow \real^n$ to be the limiting opinion vector of an SBC or SBI system whose initial opinion vector is $y$ and the interconnection topology of its agents remains unchanged for all $t \ge 0$. That is,
\begin{equation*}
  \yf(y) = \lim_{t\to\infty}A(y)^ty \in\real^n.
\end{equation*}
\end{definition}
The final value at constant topology of any equilibrium opinion vector is equal to itself, that is, $\yf(y_0)= \lim_{t\to\infty}A(y_0)^t y_0 =y_0 $. Therefore, the set of final values at constant topology is a superset of the equilibria of the system, and the set of equilibria is a superset of the limiting opinion vectors.  The condition under which a final value at constant topology is an equilibrium is discussed as follows.
\begin{proposition}[Properties of the final value at constant topology]
  \label{prop:properties-final-value}
  For any opinion vector $y\in\real^n$ in an SBC or SBI system, whose adjacency matrix can be found from equation~\eqref{canonical}:
  \begin{enumerate}
  \item $\yf(y)$ is well defined, and is equal to
    \begin{equation*}%\label{fvalueeq}
      \yf(y) = P^T(y) 
      \ml C & 0 & 0\\ 
      0 & M^*& 0 \\
      (I - \Theta)^{-1}\Theta_CC  & (I - \Theta)^{-1}\Theta_MM^* & 0\mr(y) P(y) y, 
    \end{equation*}
    where 
the submatrix $M^*(y)$ is set equal to $\lim_{t\rightarrow \infty} M(y)^t$ and is well defined.
  \item If the two networks of agents with opinion vectors $y$ and $\yf(y)$ have the same interconnection topology or, equivalently, $G_r(y) =G_r(\yf(y))$, then \label{nomoderate}
\begin{enumerate}
\item\label{iiparta} $\yf(y)$ is an equilibrium opinion vector, 
\item\label{iipartb} $G_r(y)$ contains no moderate-minded component, and
\item\label{midclosed} in any  WCC of $G_r(\yf(y))$, the maximum and minimum opinions belong to its closed-minded components.
\end{enumerate}
  \end{enumerate}
\end{proposition}
\begin{proof} 
Let us first drop the $y$ argument for matrices for readability. Regarding part~(i),
  \begin{gather*}
    \yf(y) = \lim_{t\to\infty} A^t y = P^T \lim_{t\to\infty}\overline{A}^t P y = P^T  \lim_{t\to\infty} \ml C^t & 0 & 0\\ 0 & M^t & 0 \\ \Theta_C^{(t)} & \Theta_M^{(t)} & \Theta^t \mr  P y.
  \end{gather*}
From Section~\ref{blocks}, $C^t = C$ for any $t \ge 1$, $\lim_{t\to\infty}\Theta^t = {\bf 0}$, and each diagonal block $M_i$, with size $n_i$, is a row-stochastic primitive nonnegative matrix. For such matrices the Perron-Frobenius Theorem tells us that the spectral radius is equal to one, and the essential spectral radius is strictly less than one. Thus, if we let $\nu_i \in  \real^{n_i}$ be a left eigenvector of $M_i$ for the eigenvalue one, then from \cite[Remark 1.69]{FB-JC-SM:09},
\begin{equation}\label{m_infinity} 
M_i^* = \lim_{t\to\infty} M_i^t = (\nu_i{\bf 1}_{n_i})^{-1}{\bf 1}_{n_i} \nu_i.
 \end{equation}
Using the solution to the infinite products of transition matrices of a Markov chain, given in \cite[Chapter 5]{FG:08}, it can be shown that
$\lim_{t\to\infty} \Theta_C^{(t)} = (I - \Theta)^{-1}\Theta_CC$, and $\lim_{t\to\infty}  \Theta_M^{(t)}= (I - \Theta)^{-1}\Theta_MM^*.$

Regarding part~\ref{iiparta}, $G_r(y) =G_r(\yf(y))$ results in $A(y) =A(\yf(y))$, and hence
\begin{equation*}
  A(\yf(y)) \yf(y) = A(y)  \lim_{t\to\infty}A(y)^t y =\lim_{t\to\infty}A(y)^t y = \yf(y).    
\end{equation*}

Regarding part~\ref{iipartb}, by contradiction assume that $G_r(y)$ contains at least one moderate-minded component with the opinion vector $y_{M_1}$ and the adjacency matrix $M_1$. The trajectory of each sink in $C(G_r(y))$ is independent of other nodes, hence $\yf_{M_1}(y) = \lim_{t\to\infty}M_1^t y_{M_1},$ and by equation~\eqref{m_infinity}, $\yf_{M_1}(y) =  (\nu_1{\bf 1}_{n_1})^{-1}{\bf 1}_{n_1} \nu_1 y_{M_1}$. Since $(\nu_1{\bf 1}_{n_1})^{-1}\nu_1 y_{M_1} $ is a scalar, all agents in one moderate-minded component are in consensus in final value at constant topology, and their adjacency matrix is no longer $M_1$, but rather a complete consensus matrix, which contradicts the assumption of $G_r(y) = G_r(\yf(y))$.

Regarding part~\ref{midclosed}, let $i$ denote the agent with the minimum final value at constant topology in one WCC of $G_r(\yf(y))$. By contradiction, assume that $i$ is open-minded. Granted that the confidence or influence bounds are strictly greater than zero, in the set of out-neighbors of each open-minded agent there exists at least one agent with distinct opinion. Since $i$ has the smallest opinion among its neighbors, $\yf_i(y)$ increases after taking an average of $i$'s out-neighbors opinions. In other words, the $i$th entry in the vector $A(\yf(y)) \yf(y) $ is strictly larger than $\yf_i(y)$, which contradicts the fact that $\yf(y)$ is an equilibrium opinion vector and invariant under matrix $A(\yf(y))$. Same can be proved for the agent with the maximum opinion.
\end{proof}

\section{Convergence Analysis}\label{sufficientcond}
In this section, motivated by Conjectures~\ref{xinfconj} and \ref{conjtau}, we drive sufficient condition which guarantees that an SBC or SBI trajectory converges to a limiting opinion vector. 
Next, to explain Conjecture~\ref{inffiniteconj}, we study sufficient conditions for SBC and SBI systems separately that guarantee  reaching a fixed state in finite time.

\subsection{Sufficient Condition for Constant Topology and Convergence}
Our sufficient condition is based on specific neighborhoods of each opinion vector, which is introduced in the following.

  \begin{definition}[Equi-topology distances and neighborhoods]
  Consider an SBC or SBI system with opinion vector $z\in \real^n$.
  \begin{enumerate}
  \item The \emph{equi-topology distance} of $z$ is a vector of
    non-negative entries $\epsilon(z)\in\real^n_{\geq0}$ defined by, for
    $i\in\until{n}$,
  \begin{equation}\label{epsiloneq}
    \epsilon_i(z) = 0.5 \min \{||z_i - z_j| - R| : j\in\until{n}\setminus\{i\}, R \in \{r_i,r_j\} \},
  \end{equation}
  and the \emph{equi-topology neighborhood} of $z$ is the set
  $\subscr{\mathcal{B}}{et}(z)$ of opinion vectors $y\in\real^n$ such that
  \begin{align*}
    |y_i - z_i| <  \epsilon_i(z), & \text{ for all } i\in\until{n} \text{ with }
    \epsilon_i(z) > 0, \text{ and}\\ 
    |y_i - z_i| =  \epsilon_i(z), & \text{ for all } i\in\until{n} \text{ with }
    \epsilon_i(z) = 0.
  \end{align*}

\item The \emph{invariant equi-topology distance} of $z$ is the vector of
  non-negative entries $\delta(z)\in\real^n_{\geq0}$ defined by, for
  $i\in\until{n}$,
  \begin{equation}\label{deltaeq}  
    \delta_i(z) =
    \min \{\epsilon_j(z) : j \mbox{ is a predecessor of } i \text{
        in the graph } G_r(z) \},
  \end{equation}
  and the \emph{invariant equi-topology neighborhood} of $z$ is the set
  $\subscr{\mathcal{B}}{iet}(z)$ of opinion vectors $y\in\real^n$ such
  that
  \begin{align*}
    |y_i - z_i| <  \delta_i(z), & \text{ for all } i\in\until{n} \text{ with }
    \delta_i(z) > 0, \text{ and}\\ 
    |y_i - z_i| =  \delta_i(z), & \text{ for all } i\in\until{n} \text{ with }
    \delta_i(z) = 0.
  \end{align*}
  \end{enumerate}
\end{definition}

\begin{figure}[h!]
  \begin{minipage}[b]{.33\linewidth}
    \includegraphics[width=.9\textwidth,keepaspectratio]{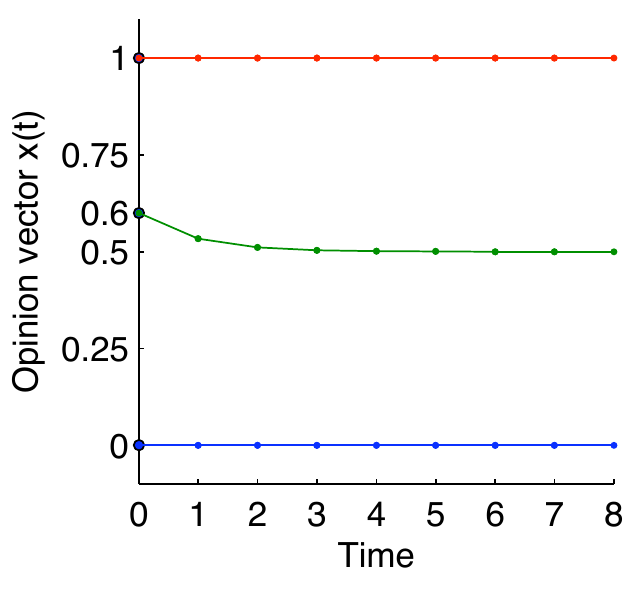}    
  \end{minipage}\hfill\begin{minipage}[b]{.67\linewidth}
    \caption{The SBC trajectory $x(t)$ with $x(0)=[0\;\; 0.6 \;\; 1]^T$ and confidence bounds $r=[0.25\;\;
      1 \;\; 0.25]^T$ is illustrated. It can be computed that $\yf(x(0)) =
      [0\;\; 0.5 \;\; 1]$ and $\delta(\yf(x(0))) = [0.25 \;\; 0.25 \;\;
      0.25]$.  Clearly, $x(0) \in \subscr{\mathcal{B}}{iet}(\yf(x(0)))$,
      i.e., the initial vector satisfies the special case of
      Theorem~\ref{epsilonthm} stated in Remark~\ref{specialcase}. Hence,
      $x(t)$ converges to $\yf(x(0))$, and their proximity digraphs are
      equal. However, if the confidence bounds are equal to $r=[0.5\;\; 1
      \;\; 0.25]^T$, then $\delta(\yf(x(0))) = [0\;\;0\;\;0]$, and $x(t)
      \notin \subscr{\mathcal{B}}{iet}(\yf(x(0)))$ for all $t\ge 0$.
      Therefore, $x(t)$ converges to $\yf(x(0))$, while their proximity
      digraphs are different. Both trajectories with the two
      confidence bounds vectors are the
      same.   \label{znotxstarfig}} 
%This example illustrates trajectories for which an equilibrium vector $z$ that satisfies the condition of the  Theorem~\ref{epsilonthm} does not exist, however $x(t)$ converges to $\yf(x(0))$ as time goes to infinity.
  \end{minipage}
\end{figure}

 Note that in any SBC or SBI digraph, each node has a self-loop, and hence each agent is a predecessor of
itself. Therefore, for any opinion vector $z\in\real^n$ and for all $i\in\until{n}$, we have $\delta_i(z) \le \epsilon_i(z)$, which results in $\subscr{\mathcal{B}}{iet}(z) \subset \subscr{\mathcal{B}}{et}(z)$.
\begin{lemma}[Sufficient condition for equal topologies]\label{epsilonlemma}
  Consider an SBC or SBI system with opinion vectors $y,z\in \real^n$. If
  $y$ belongs to the equi-topology neighborhood of $z$, then the two
  networks of agents with opinion vectors $y$ and $z$ have the same
  interconnection topology, or equivalently $G_r(y) = G_r(z)$.
\end{lemma}
\begin{remark} 
 For any $y\in\real^n$, if $y \in \subscr{\mathcal{B}}{et}(\yf(y))$, then by Lemma~\ref{epsilonlemma} we have
 $G_r(y) = G_r(\yf(y))$.
  Hence, by Proposition~\ref{prop:properties-final-value}, $\yf(y)$
  is an equilibrium opinion vector, and $G_r(y)$ contains no moderate-minded
  component.
\end{remark}
\begin{proof}
  For any $i,j\in \until{n}$ two cases exists:\\
  1. $j$ is an out-neighbor of $i$ in $G_r(z)$, hence $|z_i - z_j| \le r$,
  where in an SBC system $r=r_i$ and in an SBI system $r=r_j$. In either
  system, since $y \in \subscr{\mathcal{B}}{et}(z)$,
\begin{equation*}
|y_i-y_j| \le |z_i - z_j| +\epsilon_i(z)+\epsilon_j(z) \le |z_i - z_j|+ ||z_i - z_j|| -r| = r.
\end{equation*}
2. $j$ is not an out-neighbor of $i$ in $G_r(z)$, hence
$|z_i - z_j|> r$, with $r$ defined above. If both $\epsilon_i(z)$ and $\epsilon_j(z)$ are zero, then $y \in \subscr{\mathcal{B}}{et}(z)$ gives us
$$|y_i-y_j| = |z_i - z_j| > r,$$ 
and if at least one is nonzero, then
\begin{equation*}
|y_i-y_j| > |z_i - z_j| -\epsilon_i(z)-\epsilon_j(z) \ge |z_i - z_j| - ||z_i - z_j| -r| = r.
\end{equation*}
Therefore, the neighboring relation of agents in $G_r(z)$ is preserved in $G_r(y)$. 
One can also prove that any neighboring relation in $G_r(y)$ is preserved in $G_r(z)$.
\end{proof}
 \begin{theorem}[Sufficient condition for constant topology and convergence]\label{epsilonthm}
  Consider a trajectory $x(t)$ of an SBC or SBI system. Assume that there exists an equilibrium opinion vector $z\in\real^n$ for the system such that $x(0) \in \real^n$ belongs to the invariant equi-topology neighborhood of $z$. Then, for all $t \ge 0$:
\begin{enumerate}
\item\label{epsst} $x(t)$ takes value in the equi-topology neighborhood of
  $z$, and hence $G_r(z) = G_r(x(t))$;
\item\label{nomoderatest} $G_r(x(t))$ contains no moderate-minded component; and
\item\label{convergencest} $x(t)$ converges to $\yf(x(0))$ as time goes to infinity. %$t\rightarrow \infty$.
\end{enumerate}
\end{theorem}
\begin{remark}[Interpretation of Theorem~\ref{epsilonthm}]
This theorem tells us that if the trajectory of an SBC or SBI system enters a specific ball around any equilibrium opinion vector of that system, then it remains in some  larger ball around that vector for all future iterations. Moreover, the proximity digraph of the trajectory and the equilibrium opinion vector remain equal.
\end{remark}
\begin{remark}\label{remarknotgeneral}
Under the condition of Theorem~\ref{epsilonthm}, a trajectory $x(t)$ converges to its final value at constant topology $\yf(x(t))=\yf(x(0))$. However, $\yf(x(t))$ is not necessarily equal to the equilibrium opinion vector $z$, and the proximity digraphs $G_r(\yf(x(t)))$ and $G_r(x(t))$ can be different, see Figure~\ref{znotxstarfig}.
 \end{remark}
\begin{remark}\label{specialcase}
One special case of Theorem~\ref{epsilonthm} is when $x(0) \in \subscr{\mathcal{B}}{iet}(\yf(x(0))) $, which implies that $\yf(x(0))$ is an equilibrium opinion vector, again see Figure~\ref{znotxstarfig}.
\end{remark}
\begin{proof}[Proof of Theorem~\ref{epsilonthm}]
Regarding statement~\ref{epsst}, by induction we prove that $x(t) \in \subscr{\mathcal{B}}{et} (z)$ for all $t \ge 0$, which by Lemma~\ref{epsilonlemma} results in $G_r(x(t)) = G_r(z)$. 
The first induction step is $x(0) \in \subscr{\mathcal{B}}{et} (z)$, which is true knowing that $x(0) \in \subscr{\mathcal{B}}{iet} (z)$ and $\subscr{\mathcal{B}}{iet} (z) \subset \subscr{\mathcal{B}}{et} (z)$. 
To complete the induction argument, assume that the statement~\ref{epsst} holds at times $t =0 ,\dots, \tau$, which implies that $A(z) = A(x(t))$.
The equilibrium opinion vector $z$ satisfies $z=A(z)z$, thus we have $x(t+1) - z = A(x(t))x(t) - A(z)z = A(z)(x(t) - z)$ or, equivalently
\begin{equation}\label{deltadistancekids2}
x_i(t+1) - z_i = \frac{1}{|\mathcal{N}_i(z)|} \sum_{j\in \mathcal{N}_i(z)} (x_j(t) - z_j).
\end{equation}
One can see that 
\begin{multline}\label{ineqkids2}
|x_i(\tau + 1) - z_i| \le \max_{j \in \mathcal{N}_i(z)} |x_j(\tau) - z_j| \le \max_{\ell \in \mathcal{N}_j(z), j \in \mathcal{N}_i(z)} |x_\ell(\tau-1) - z_\ell|\\  
\le  \dots \le \max_{k\in \mathcal{M}} |x_k(0) - z_k|,
\end{multline}
where $\mathcal{M}$ is a subset of successors of $i$ in $G_r(z)$, and thus for any $k\in\mathcal{M}$, equation~\eqref{deltaeq} tells us that $\delta_k(z) \le \epsilon_i(z)$.
Here again two cases exists: First, if for all $k\in\mathcal{M}$, $\delta_k(z) = 0$, then the condition $x(0) \in \subscr{\mathcal{B}}{iet} (z)$ implies that $x_k(0) - z_k = 0$, and it follows from inequality~\eqref{ineqkids2} that $x_i(\tau+1) - z_i =0$. Second, if there exists $\ell \in \mathcal{M}$ such that $\delta_\ell(z) > 0$, then $\epsilon_i(z) > 0$ and  
$$|x_i(\tau+1) - z_i| \le  \max_{k\in \mathcal{M}} |x_k(0) - z_k| < \max_{k\in \mathcal{M}} \delta_k(z) \le \epsilon_i(z).$$ 
Therefore, $x(\tau+1) \in \subscr{\mathcal{B}}{et} (z)$.

Regarding statement~\ref{nomoderatest}, according to Section~\ref{final_value}, an equilibrium opinion vector is equal to its own final value at constant topology. Hence $G_r(z) = G_r(z^*(z))$, and by Proposition~\ref{prop:properties-final-value}, $G_r(z)$ and thus $G_r(x(t))$ contain no moderate-minded component. 

Regarding statement~\ref{convergencest}, according to the definition of the final value at constant topology, if the topology remains constant for all $t \ge 0$, then $x(t)$ converges to $\yf(x(0))$. 
\end{proof}

Motivated by Conjecture~\ref{xinfconj}, the existence of a limiting opinion vector is required in the following lemma.
\begin{lemma}[Sufficient condition for a limiting opinion vector to be an equilibrium]\label{xinflemma}
Pick a trajectory $x(t)$ of an SBC or SBI system that is convergent. Denote the limiting opinion vector of $x(t)$ by $x_\infty$. If $\min_{i \in \until{n}} \epsilon_i(x_\infty) >0$, where $\epsilon(x_\infty)$ is the equi-topology distance of $x_\infty$ , then there exists time $T$ such that for all $t\ge T$:
\begin{enumerate}
\item $G_r(x_\infty) = G_r(x(t))$, and 
%\item $\|x(t) - \yf(x(T))\|_\infty < \min_{i \in \until{n}}\epsilon_i(\yf(x(T)))$ for all $t\ge T$, and 
\item $x_\infty = \yf(x(t))$, and is an equilibrium opinion vector.
\end{enumerate}
\end{lemma}
\begin{proof}
According to the definition of convergence, for any $\delta \in \real_{>0}$, there exists $T$ such that for all $t \ge T$, $\|x(t) - x_\infty \|_\infty < \delta$.
Now, if we let $\delta$ be equal to $\min_{i \in \until{n}} \epsilon_i(x_\infty)$, then $\|x(t) - x_\infty \|_\infty <  \min_{i \in \until{n}} \epsilon_i(x_\infty)$ for all $t \ge T$, and it follows from Lemma~\ref{epsilonlemma} that $ G_r(x_\infty) = G_r(x(t))$.  Under fixed topology, $x(t)$ converges to its final value at constant topology, thus $x_\infty = \yf(x(t))$. Moreover, the equality $G_r(x(t)) = G_r(\yf(x(t))) $ tells us that $\yf(x(t))$, and hence $x_\infty$, is an equilibrium opinion vector.
%Based on the definition of convergence, for any $\delta > 0$, there exists time $T$ such that $\|x(t) - x_\infty \|_\infty < \delta$ for all $t \ge T$. If we let $\delta = \min_{i \in \until{n}} \epsilon_i(x_\infty)$, then 
\end{proof}
\begin{corollary} %[Lyapunov stable equilibria]
An equilibrium opinion vector $z$ is a Lyapunov stable equilibrium vector for the system if  $\min_{i \in \until{n}} \epsilon_i(z) > 0$.
\end{corollary}

\subsection{Sufficient Condition for Convergence in Finite Time and Consensus}\label{consensus}

In this subsection, we discuss the sufficient conditions for SBC and SBI
systems to converge to an \emph{agreement opinion vector}. In an agreement
opinion vector, any two agents are either disconnected or in consensus. Reaching global consensus, in which all agents hold the same opinion, is an special case of convergence to an agreement opinion vector. Note that it is possible for an SBC or SBI system to reach an opinion vector that
contains two neighbor agents with separate opinions in finite time. For instance, the trajectory of an SBC system with initial opinion vector  $ [0  \;\;    2  \;\;   3 \;\;    4.5\;\;  7 ]^T$ and bounds vector $  [0.01    \;\;   3  \;\;   0.01  \;\;  3  \;\;   0.01   ]^T$ exhibits convergence to a fixed profile in finite time. While, the SBC digraph of the limiting opinion vector contains open-minded agents.  
\begin{proposition}[Properties of agreement opinion vectors]
For any agreement opinion vector $\tilde{y} \in \real^n$ in an SBC or SBI system: 
\begin{enumerate}
\item\label{epsilongezero} $\min_{i \in \until{n}} \epsilon_i(\tilde{y}) >0$, where $\epsilon(\tilde{y})$ is the equi-topology distance of $\tilde{y}$; and
\item\label{agreementfinite} if $\tilde{y}$ is the limiting opinion vector of a trajectory, then the trajectory reaches $\tilde{y}$ in finite time.
\end{enumerate}
\end{proposition}
\begin{proof}
  Regarding statement~\ref{epsilongezero}, by contradiction assume that $\min_{i
    \in \until{n}} \epsilon_i(\tilde{y}) =0$. Then based on
  equation~\eqref{epsiloneq}, there exist agents $i$ and $j$ such that
  $|\tilde{y}_i-\tilde{y} _j| = r_i$.  The latter equation tells us that $j
  \in \mathcal{N}_i(\tilde{y})$ in an SBC digraph or $i \in
  \mathcal{N}_j(\tilde{y})$ in an SBI digraph, while their opinions are
  different from each other by $r_i$, which contradicts the definition of
  agreement opinion vectors. Regarding statement~\ref{agreementfinite},
 consider trajectory $x(t)$ that converges to $\tilde{y}$. Then,
  previous statement shows that the limiting opinion vector of $x(t)$ satisfies the condition of
  Lemma~\ref{xinflemma}. Therefore, there exists time step $\tau$ such that
  the proximity digraphs $G_r(\tilde{y})$ and $G_r(x(t))$ are equal for all
  $t \ge \tau$. On the other hand, the proximity digraph of an agreement
  opinion vector contains only closed-minded components. Hence, the agents
  in each WCC of $G_r(x(\tau))$ reach consensus at the next iteration.
\end{proof}

One sufficient condition that guarantees asymptotic consensus in
``agreement algorithms'', which includes SBC and SBI systems, is given in \cite[Theorem 2.4]{AO-JNT:07} and is as follows. Take a trajectory $x(t)$ of an SBC or SBI
system with proximity digraph $G_r(x(t)) = \big(V, E(x(t))\big)$, where $V$ and
$E(x(t))$ are the sets of nodes and edges of the digraph,
respectively. If there exists $\tau$ such that the graph $\big( V,
E(x(k\tau))\cup E(x(k\tau+1))\cup \dots \cup E(x((k+1)\tau-1))\big)$ is
strongly connected for all $k \in \mathbb{Z}_{\ge 0}$, then all entries of $x(t)$ converge to one real number.
However, this sufficient condition requires knowledge of the system for
infinite time, and is the same for both SBC and SBI systems.  Hence, we
derive sufficient conditions that are required to hold in one time step, and
also make it possible to compare SBC and SBI systems in support of
Conjecture~\ref{inffiniteconj}. Let us first define the \emph{opinion
  interval} of any subgraph of an SBC or SBI digraph be a closed
interval in $\real$ between that subgraph's minimum and maximum opinions.

\begin{proposition}[Sufficient conditions for convergence to an agreement opinion vector]\label{agreementprop}
Consider the opinion vector $y\in\real^n$ in an SBC or SBI system with the following properties:
\begin{enumerate}
\item\label{separateintervals} the opinion intervals of any two WCC's of the proximity digraph are separated from each other by a distance strictly larger than the maximum confidence or influence bounds of the agents in those WCC's; and
\item\label{sinkforever} it is true that:
\begin{itemize}
\item for any WCC of $y$'s SBC digraph, with $m$ agents, at least $m-1$ agents have confidence bounds larger than that WCC's opinion interval; and
\item for any WCC of $y$'s SBI digraph, at least one agent has influence bound larger than that WCC's opinion interval.
\end{itemize}
\end{enumerate}
Then, the trajectories of both SBC and SBI systems with the initial opinion vector $y$ converge to agreement opinion vectors in finite time. Moreover, in every WCC of either of the SBC or SBI digraphs, at least one node is an out-neighbor of all nodes in that WCC for all $t \ge 0$. 
\end{proposition}
\begin{remark}
Any trajectory of an SBC or SBI system that converges to an agreement opinion vector will eventually satisfy the conditions of Proposition~\ref{agreementprop}.
\end{remark}
\begin{proof}[Proof of Proposition~\ref{agreementprop}] 
 Let us denote either of the SBC or SBI digraphs of $y$ by $G_r(y)$.
In an SBC or SBI system, the smallest and largest opinions in a separate WCC of the proximity digraph are, respectively, non-decreasing and non-increasing in one iteration \cite{VDB-JMH-JNT:09}. This fact tells us that for all $t\ge0$:
first, under the condition~\ref{separateintervals}, the two sets of nodes of two separate WCC's in $G_r(x(0))$ remain separate in $G_r(x(t))$;
second, if the condition~\ref{sinkforever} holds for $G_r(x(0))$, then it also holds for $G_r(x(t))$.
Now, under condition~\ref{sinkforever} for both SBC and SBI systems, any WCC of $G_r(x(t))$ contains at least one agent that is an out-neighbor of all agents in that WCC for all $t \ge 0$. 
%This topology results in a strict decrease in the opinion interval of each WCC in one iterations. 
Denote one such agent in a WCC by $s$, then that WCC's agents with maximum and minimum opinions update their opinions by taking an average of their out-neighbors, including $s$. Hence, at the next iteration, their opinions will converge to $s$'s opinion, which results in an strict decrease in the opinion interval of the WCC. Since the confidence or influence bounds are strictly greater than zero, there exists a time step after which the opinion interval of the WCC is larger than the minimum confidence or influence bound. Consequently, all agents become each others out-neighbors, and the WCC becomes one closed-minded component.
\end{proof}

\section{Numerical Analysis}\label{expsec}

\begin{figure}
  \centering
  \includegraphics[width=.49\textwidth,keepaspectratio]{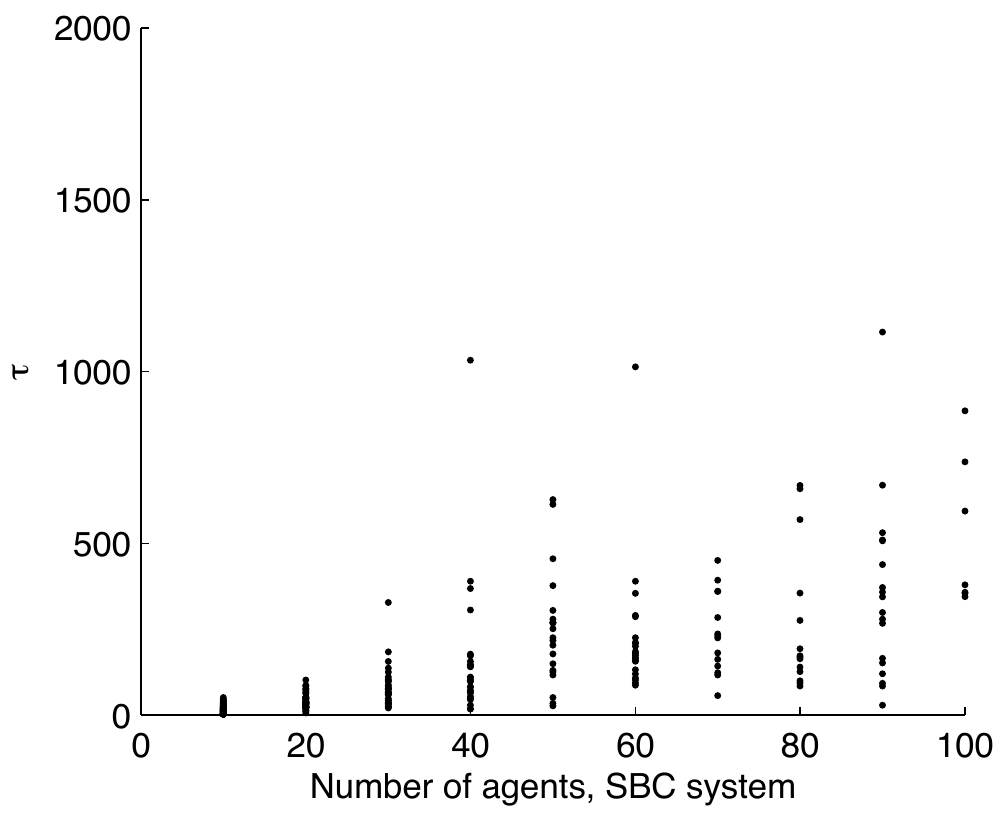} \hfill
  \includegraphics[width=.49\textwidth,keepaspectratio]{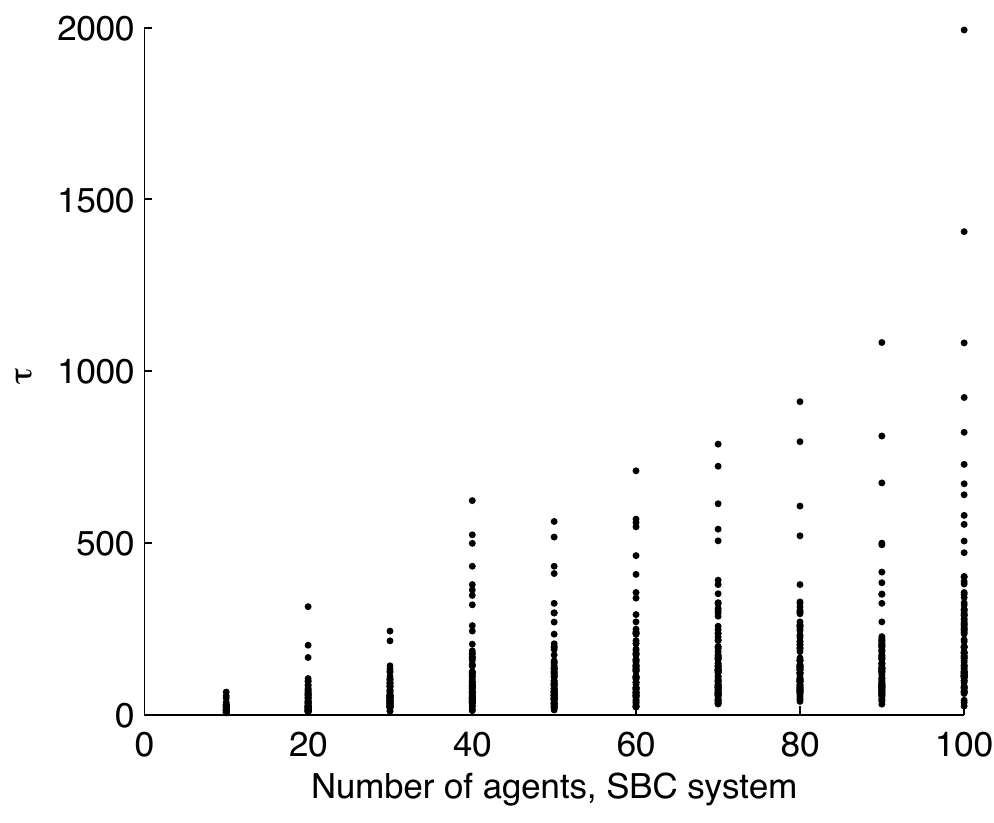}
\caption{In one thousand simulations of SBC systems, the time step $\tau$ at which the trajectory of each system satisfied the sufficient condition of Theorem~\ref{epsilonthm} is plotted versus the number of agents in that system. The left and right plots, respectively, illustrate time $\tau$ for trajectories that converged in finite time and infinite time.
Each initial opinion vector and bounds vector are generated randomly and uniformly distributed on [0, 1] and [0, 0.3], respectively. 
For each agent number hundred simulations are performed. All trajectories satisfied the special case of sufficient condition of Theorem~\ref{epsilonthm}, stated in Remark~\ref{specialcase} in finite time.}\label{tfConf}
\end{figure}

\begin{figure}
  \centering
  \includegraphics[width=.49\textwidth,keepaspectratio]{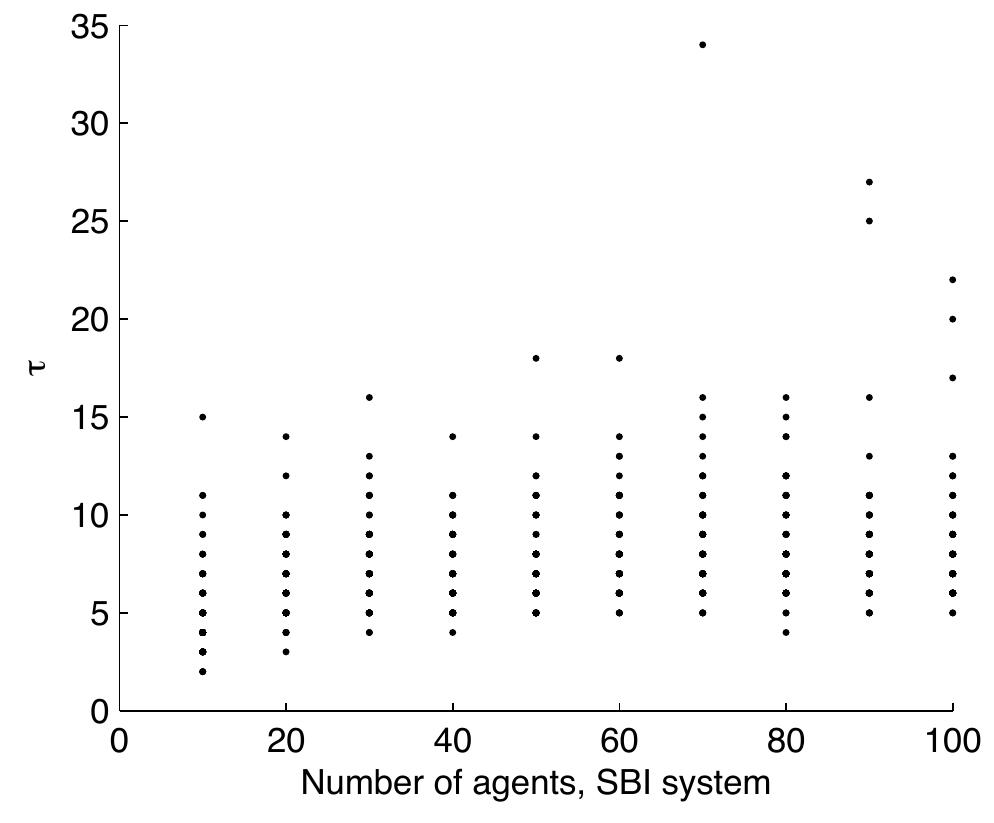} \hfill
  \includegraphics[width=.49\textwidth,keepaspectratio]{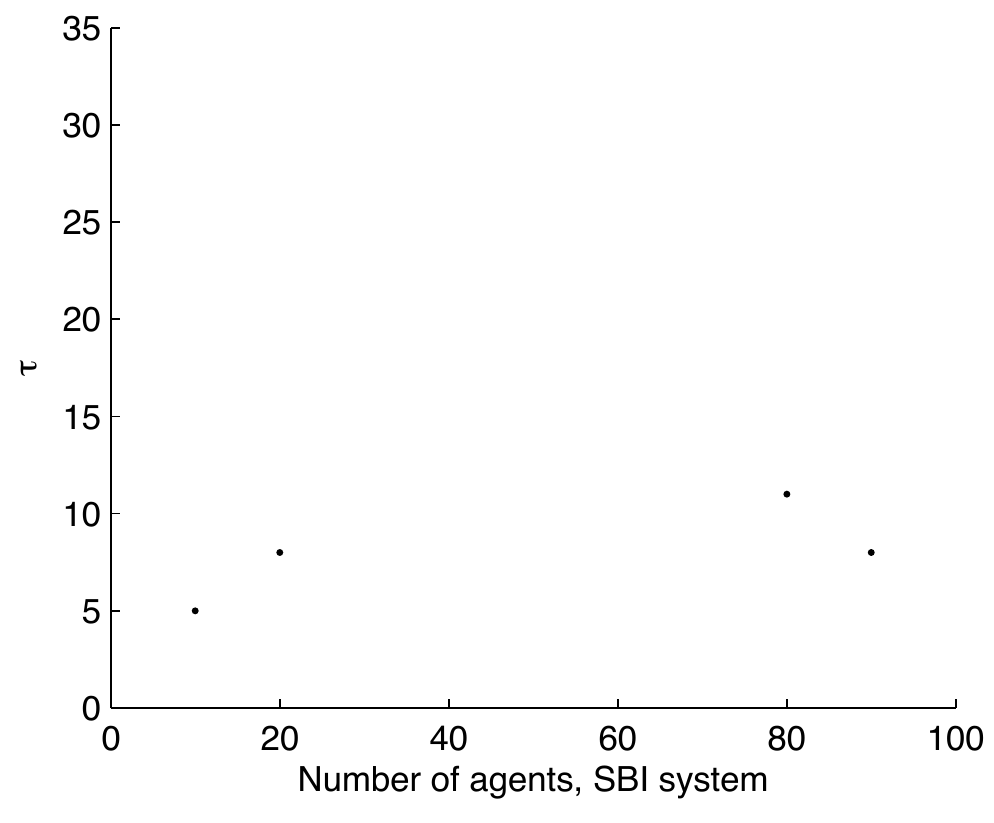}
  \caption{In one thousand simulations of SBI systems, the time step $\tau$ at which the trajectory of each system satisfied the sufficient condition of Theorem~\ref{epsilonthm} is plotted versus the number of agents in that system. The left and right plots, respectively, illustrate time $\tau$ for trajectories that converged in finite time and infinite time. As shown, only four SBI trajectories converged in infinite time.
Each initial opinion vector and bounds vector are generated randomly and uniformly distributed on [0, 1] and [0, 0.3], respectively. 
For each agent number hundred simulations are performed. All trajectories satisfied the special case of sufficient condition of Theorem~\ref{epsilonthm}, stated in Remark~\ref{specialcase} in finite time.}\label{tfInfl}
\end{figure}

In this section, we provide extensive simulation results that demonstrate the results of Section~\ref{sufficientcond} and are consistent with our conjectures.
 We performed 2000 simulations: 100 simulations of both SBC and SBI systems for
ten different agent numbers. In each
simulation, the initial opinion vector and 
bounds vector are generated randomly and uniformly distributed on
[0, 1] and [0, 0.3], respectively.
The time steps $\tau$ at which trajectories satisfied the condition of Theorem~\ref{epsilonthm} are plotted 
in Figures~\ref{tfConf} and \ref{tfInfl}.
All the 2000 SBC and SBI trajectories eventually satisfied the special case of the
sufficient condition of Theorem~\ref{epsilonthm}, stated in
Remark~\ref{specialcase}.  In other words, for each trajectory $x(t)$,
there exists time $\tau$ such that $x(\tau)$ belongs to the invariant
equi-topology neighborhood of its own final value at constant topology
$\yf(x(\tau))$. Thus, $\yf(x(\tau))$ is an equilibrium opinion vector, and is
equal to the limiting opinion vector of $x(t)$. The frequency of occurrence
of this special case is intuitively explained by the following statements: First, by
Conjecture~\ref{xinfconj}, for each trajectory a limiting opinion vector
$x_\infty$ exists. Second, for any randomly generated opinion vector $y$
and bounds vector $r$, the probability of having
$\min_{i\in \until{n}} \epsilon_i(y) = 0$, where $\epsilon(y)$ is the
equi-topology distance of $y$, is equal to zero, and one can assume that
the same holds for any limiting opinion vector. Third,
Lemma~\ref{xinflemma} tells us that if the limiting opinion vector
satisfies $\min_{i\in \until{n}} \epsilon_i(x_\infty) > 0$, then the trajectory eventually satisfies the mentioned special case of condition of Theorem~\ref{epsilonthm}.

\begin{figure}
  \centering
\includegraphics[width=.57\textwidth,keepaspectratio]{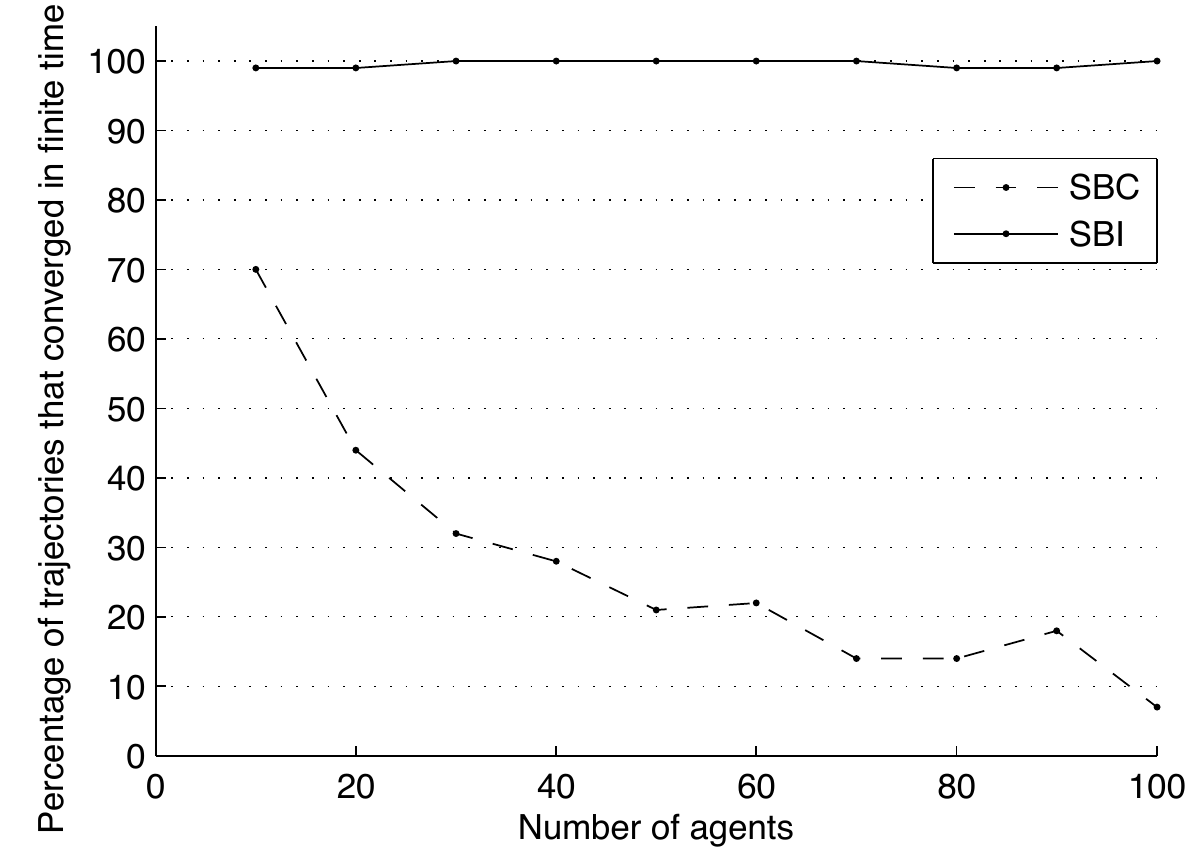}
  \caption{For each agent number 100 SBC systems and 100 SBI systems are simulated, as explained in Figures~\ref{tfConf} and \ref{tfInfl}.
The percentage of SBC (blue) and SBI (green) trajectories that reached their limiting opinion vector, which is an agreement opinion vector, in finite time is plotted versus agent number.}\label{percentfig}
\end{figure}
In above mentioned simulations, for each agent number, the percentage of SBC and SBI trajectories that reached a fixed profile in finite time are plotted in Figure~\ref{percentfig}. Clearly, Figure~\ref{percentfig} supports Conjecture~\ref{inffiniteconj}. To explain this frequency of convergence of SBI trajectories in finite time as compared with SBC trajectories, we use the results of Subsection~\ref{consensus}. 
For uniformly randomly generated opinion vector and bounds vector, an SBI digraph is more likely to satisfy condition~\ref{separateintervals} of Proposition~\ref{agreementprop} than an SBC digraph. 
One can assume that the same holds for a trajectory with uniformly randomly generated initial opinion vector and bounds vector.

In the next section, based on our ``constant topology in finite time'' conjecture, we assume that the interconnection topology in an SBC or SBI system remains constant for infinite time and address the following questions: How the three classes of agents behave? How groups of agents affect each other? And can one explain the observed pseudo-stable behavior of trajectories, as stated in Conjecture~\ref{conjtau}?

\section{The Rate and Direction of Convergence under Fixed Topology}\label{fixedtopbehavior}

In this section, we analyze the rates and directions of convergence of separate classes of agents in the SBC and SBI systems under fixed interconnection topology as time goes to infinity. This analysis proves that the system shows a pseudo-stable behavior under fixed topology.
\begin{definition}[Agent's per-step convergence factor]\label{k}
In an SBC or SBI system with trajectory $x(t) \in\real^n$, we define the \emph{per-step convergence factor} of an agent $i$ whose $x_i(t) - \yf_i(x(t))$ is nonzero to be 
\begin{equation*}
k_i(x(t)) = \frac{x_i(t+1) - \yf_i(x(t)) }{ x_i(t) - \yf_i(x(t))}. 
\end{equation*}
\end{definition}

The per-step convergence factor of a network of agents was previously introduced in~\cite{LX-SB:04} to measure the overall speed of convergence toward consensus.% in order to solve the problem of finding a linear iteration that yields distributed averaging over the network.
\begin{remark}\label{pseudoperstep}
Consider a converging trajectory $x(t)$ whose limiting opinion vector is equal to $\yf(x(t))$ for all $t \ge 0$. Then, $x(t)$ exhibits a pseudo-stable behavior, see equation~\eqref{pseudoeq}, if and only if for all $i \in \until{n}$
\begin{equation*}
\begin{cases}
0 < k_i(x(t)) < 1,  \hspace{.6in}  & \mbox{if} \hspace{.1in}  k_i(x(t)) \mbox{ exists}, \\
x_i(t)  =  x_i(t+1) = \yf_i(x(t)), \hspace{.6in} &   \mbox{otherwise}.
\end{cases}
\end{equation*}
\end{remark}
\begin{definition}[Leader SCC]\label{leaderSCC}
For opinion vector $y\in\real^n$, let $G_r(y)$ denote either its SBC or SBI digraph.
For any open-minded SCC of $G_r(y)$, $S_k(y)$, denote the set of its open-minded successor SCC's by $\mathcal{M}(S_k(y))$, which includes $S_k(y)$.
We define $S_k(y)$'s \emph{leader SCC} to be the SCC whose adjacency matrix has the largest spectral radius among all SCC's of $\mathcal{M}(S_k(y))$. 
\end{definition}

Note that in SBC and SBI digraphs, the adjacency matrix of a large SCC has a large spectral radius, hence that SCC tends to become a leader SCC for its predecessors.
\begin{theorem}[Evolution under constant topology]\label{convlemma}
Consider an SBC or SBI system, denote its trajectory by $x(t)$ and proximity digraph by $G_r(x(t))$.
Assume that there
  exists a time $\tau$ after which $G_r(x(t))$ remains
  unchanged, that is, $G_r(x(t)) = G_r(x(\tau))$. Then, the following statements hold: 
  \begin{enumerate}
  \item\label{convfirst} $\yf(x(t)) = \yf(x(\tau))$ for all $t \ge \tau$.
  \item\label{convsecond} $G_r(x(\tau))$ contains no moderate-minded component.
  \item\label{convthird} Consider any open-minded SCC $S_k(x(t))$ of $G_r(x(t))$ and its leader SCC $S_m(x(t))$, with adjacency matrices denoted by $\Theta_{k}$ and $\Theta_{m}$, respectively. Then,
\begin{enumerate}
\item\label{convthirda} for any $i \in S_k(x(t))$, either $x_i(t) - \yf_i(x(t)) = 0$ or its per-step convergence factor converges to the spectral radius of $\Theta_{m}$ as time goes to infinity, and
\item\label{convthirdb} if the spectral radius of $\Theta_{k}$ is strictly less than that of $\Theta_{m}$, then there exists $t_1 \ge \tau$ such that for all $i \in S_k(x(t))$, $j \in S_m(x(t))$, and $t \ge t_1$,
\begin{gather*}
      x_j(t_1) < \yf_j(x(t_1)) \quad\implies\quad  x_i(t) \leq \yf_i(x(t)),
      \\
      x_j(t_1) > \yf_j(x(t_1)) \quad\implies\quad  x_i(t) \geq \yf_i(x(t)).
    \end{gather*} 
 \end{enumerate}
\item\label{convforth} There exists time $t_2 \ge \tau$ such that for all $t \ge t_2$, $x(t)$ exhibits a pseudo-stable behavior, see equation~\eqref{pseudoeq}.
\end{enumerate}
\end{theorem}
\begin{remark}[Interpretation of statement~\ref{convthird} in Theorem~\ref{convlemma}]
Parts (a) and (b) tell us, respectively, that the rates and directions of convergence of opinions in an open-minded SCC toward the final value at constant topology are governed by the rate and direction of convergence of its leader SCC.
It is easy to see that the per-step convergence factor has an inverse relation with the rate of convergence to the final value at constant topology. Therefore,  Theorem~\ref{convlemma} implies that under fixed interconnection topology, individuals converge to a final decision as slow as the slowest group of agents whom they listen to. 
\end{remark}
\begin{proof}
Statement \ref{convfirst} is a direct consequence of $A(x(t)) = A(x(\tau))$ for all $t \ge \tau$. 
Statement \ref{convsecond} can be proved similar to part~\ref{iipartb} of Proposition~\ref{prop:properties-final-value}. It was shown that under fixed interconnection topology, all agents in one moderate-minded SCC of an SBC or SBI digraph reach consensus as time goes to infinity. Since, the bounds vector is strictly greater than zero, there exists a time step after which the adjacency matrix of one moderate-minded SCC transforms into a complete consensus matrix, which contradicts the assumption of having fixed topology for infinite time.
Before proving statement~\ref{convthird}, 
since the canonical permutation matrix remains unchanged, let us assume that the opinions in $x(\tau)$ are ordered such that $A(x(\tau))=\overline{A}(x(\tau))$. Furthermore, owing to the fixed interconnection topology, we drop the $x(t)$ argument for simplicity.
Therefore, by equation~\eqref{canonical}, $$A(x(t)) = \ml C & 0 \\ \Theta_C & \Theta \mr.$$
Now, for all $t > \tau$ we have
\begin{equation*}
x(t) - \yf (x(\tau))=   \ml  C x_C(\tau) \\ x_\Theta(t) \mr - \ml  C x_C(\tau) \\ \yf_\Theta(x(\tau)) \mr 
= \ml 0 \\  x_\Theta(t) - \yf_\Theta(x(\tau)) \mr, 
\end{equation*}
where $x_C(t)$ and $x_\Theta(t)$ are the opinion vectors of agents in closed- and open-minded classes respectively. Using $\yf(x(\tau)) = A(x(\tau)) \yf(x(\tau))$, the following recurrence relation holds
\begin{equation}\label{WCCdelta}
 x_\Theta(t+1) - \yf_\Theta(x(\tau)) =   \Theta ( x_\Theta(t) - \yf_\Theta(x(\tau))) \hspace{.2in}\forall~t\ge\tau. 
\end{equation}
Consider an open-minded WCC of $G_r(x(t))$, denoted by $W_1$. Let $\Theta_1$ denote $W_1$'s adjacency matrix, and 
$x_1(t)$ denote the trajectory of nodes of $W_1$.
Under fixed interconnection topology, the trajectory of each WCC is independent of others, thus for all $t \ge 0$,
$x_1(t+\tau) - \yf_1(x(\tau)) = \Theta_1^t ( x_1(\tau) - \yf_1(x(\tau))).$
According to the block lower triangular from of $\Theta_1$, %As mentioned in Section~\ref{blocks},
\begin{equation*}
\Theta_1^t= \ml \Theta_{11}^t &   & 0\\ \Theta_{21}^{(t)}  & \Theta_{22}^t  & \\ \vdots &  &   \ddots \mr,
\end{equation*}
where each $\Theta_{ii}$ is the adjacency matrix of an SCC, denoted by $S_{ii}$, of $W_1$. Let $x_{ii}(t)$ be the opinion trajectory of nodes of $S_{ii}$. Clearly, $S_{11}$ is one of the sink SCC's in $W_1$, and $\Theta_{ii}$'s are ordered in $\Theta_1$ according to the distance of $S_{ii}$'s to the sinks. For simplicity, we prove statement~\ref{convthird} for $S_{11}$ and an SCC that is the direct predecessor of $S_{11}$. 
Without loss of generality, let $S_{22}$ be one such SCC. The proof for the rest of open-minded SCC's is similar.

Each block $\Theta_{ii}$ is nonnegative and primitive. By Perron-Frobenius Theorem: the spectral radius of $\Theta_{ii}$, denoted by $\lambda_i$, is positive and a simple eigenvalue of $\Theta_{ii}$; and there exists a positive eigenvector $\nu_i$ for $\Theta_{ii}$ associated to $\lambda_i$. Any $\Theta_{ii}$  can be written in Jordan normal form by some similarity transformation
\begin{equation*}
\Theta_{ii}  =  QJQ^{-1} = \ml \nu_i \hspace{.1in}\vline   & Q_e \mr \ml \lambda_i & {\bf 0} \\ {\bf 0} & J_e \mr \ml  w_i   \\\hline  Q_e^{(-1)} \mr,
\end{equation*}
where $w_i $ is the first row of $Q^{-1}$. Consequently, 
\begin{equation}\label{SCCconv}
\lim_{t\rightarrow \infty}\Theta_{ii}^t = \lim_{t\rightarrow \infty} (\lambda_i^t \nu_i w_i+ Q_e J_e^tQ_e^{(-1)}) = \lim_{t\rightarrow \infty} \lambda_i^t \nu_i w_i.
\end{equation}

In any open-minded WCC, the sink SCC is not effected by other SCC's, hence a sink SCC is its own leader. Therefore, for all $t \ge 0$, $$x_{11}(t+\tau) - \yf_{11}(x(\tau)) = \Theta_{11}^t ( x_{11}(\tau) - \yf_{11}(x(\tau))).$$
In the interest of simplicity, let us denote the vector $x_{ii}(t)- \yf_{ii}(x(\tau))$ by $\Delta_{ii}(t)$, then we have
\begin{equation}\label{theta11eq}
\lim_{t \rightarrow \infty} \Delta_{11}(t) =  \nu_1 w_1\Delta_{11}(\tau)\lim_{t\rightarrow \infty}  \lambda_1^t.
\end{equation}
Regarding part~\ref{convthirda} for $S_{11}$, for the per-step convergence factor of any $i \in S_{11}$ we have
\begin{equation*}
\lim_{t\rightarrow \infty} k_i(x(t)) = \lim_{t\rightarrow \infty}\frac{\lambda_1^{t+1} w_1\Delta_{11}(\tau)\nu_{1i}}{\lambda_1^t w_1\Delta_{11}(\tau)\nu_{1i}} = \lambda_1.
\end{equation*} 
Regarding part~\ref{convthirdb} for $S_{11}$, since $\lambda_1^t w_1\Delta_{11}(\tau)$ is a scalar, $\lambda_1$ is positive, and $\nu_1$ is a positive vector, all entries of vector $\nu_1 w_1\Delta_{11}(\tau) \lambda_1^t$ have the same sign. Therefore, there exists time $T \ge \tau$ after which all entries of $\Delta_{11}(t)$ have the same sign for all $t\ge T$.

Here, we prove the two statements for $S_{22}$. It can be computed that for all $t \ge 0$
\begin{gather*}
\Delta_{22}(t+\tau)=  \sum_{i=0}^{t-1} \Theta_{22}^i \Theta_{21} \Theta_{11}^{t-i-1} \Delta_{11}(\tau) + \Theta_{22}^t  \Delta_{22}(\tau).
\end{gather*}
Now, to find $\lim_{t \rightarrow \infty}\Delta_{22}(t)$, we consider three cases:\\
1) If $\lambda_1 > \lambda_2$, then $S_{11}$ is $S_{22}$'s leader. According to the transient analysis of the reducible Markov chains from \cite[Section 5.6]{FG:08}, and granted that $\lambda_1$ and $\lambda_2$ are strictly less than one, see Section~\ref{blocks},
\begin{multline*}
\lim_{t\rightarrow \infty}\sum_{i=0}^{t-1}\Theta_{22}^i\Theta_{21}\Theta_{11}^{t-i-1}=\big(\lim_{t\rightarrow \infty}\sum_{i=0}^{t-1} \Theta_{22}^i \big)\Theta_{21}\lim_{t\rightarrow \infty} \Theta_{11}^t  \\ = (I -\Theta_{22} )^{-1} \Theta_{21}\lim_{t\rightarrow \infty}\lambda_1^t \nu_1 w_1.
\end{multline*}
Therefore,
\begin{equation}\label{theta22eq}
\lim_{t \rightarrow \infty}\Delta_{22}(t) =  (I -\Theta_{22} )^{-1} \Theta_{21} \nu_1 w_1 \Delta_{11}(\tau)\lim_{t\rightarrow \infty} \lambda_1^t. 
\end{equation}
Regarding part~\ref{convthirda} for $S_{22}$, for any $i\in S_{22}$ we have
\begin{gather*}
\lim_{t \rightarrow \infty} k_i(x(t)) = \lim_{t \rightarrow \infty} \frac{\lambda_1^{t+1} w_1 \Delta_{11}(\tau) [(I -\Theta_{22} )^{-1} \Theta_{21} \nu_1]_i }{\lambda_1^t w_1 \Delta_{11}(\tau)[(I -\Theta_{22} )^{-1} \Theta_{21} \nu_1]_i} = \lambda_1.
\end{gather*}
Regarding part~\ref{convthirdb} for $S_{22}$, since $(I-\Theta_{22} )^{-1} \Theta_{21} \nu_1 $ is a nonnegative vector, all entries of the vector on the right hand side of equations \eqref{theta11eq} and \eqref{theta22eq} have the same sign as the scalar $w_1 \Delta_{11}(\tau)$.\\
2) If $\lambda_1 < \lambda_2$, then $S_{22}$ is its own leader. Similarly, 
\begin{equation}\label{theta22selfleadeq}
\lim_{t \rightarrow \infty}\Delta_{22}(t) = w_2\big( \Theta_{21} (I -\Theta_{11} )^{-1} \Delta_{11}(\tau)+\Delta_{22}(\tau)\big)\nu_2\lim_{t \rightarrow \infty}\lambda_2^t ,
\end{equation}
where $ w_2\big( \Theta_{21} (I -\Theta_{11} )^{-1} \Delta_{11}(\tau)+\Delta_{22}(\tau)\big)$ is a scalar, $\lambda_2$ is positive, and $\nu_2$ is a positive vector. Regarding part~\ref{convthirda} for $S_{22}$, for any $i \in S_{22}$ we have
\begin{gather*}
\lim_{t \rightarrow \infty} k_i(x(t)) =\lim_{t \rightarrow \infty} \frac{\lambda_2^{t+1} w_2(\Theta_{21} (I -\Theta_{11} )^{-1}  \Delta_{11}(\tau)+ \Delta_{22}(\tau)) \nu_{2i}}{\lambda_2^t w_2 (\Theta_{21} (I -\Theta_{11} )^{-1} \Delta_{11}(\tau)+\Delta_{22}(\tau) ) \nu_{2i}} = \lambda_2.
\end{gather*}
Regarding part~\ref{convthirdb} for $S_{22}$, all entries of the vector on the right hand side of equation~\eqref{theta22selfleadeq} have the same sign.\\
3) if $ \lambda_1 = \lambda_2 = \lambda$, then we have
\begin{gather*}
\lim_{t \rightarrow \infty} \Delta_{22}(t)=  \big(\alpha\nu_2 +\beta(I -\Theta_{22} )^{-1} \Theta_{21}\nu_1 \big)\lim_{t\rightarrow \infty} \lambda^t,
\end{gather*}
where $\beta=w_1 \Delta_{11}(\tau)$ and $\alpha = w_2 \Theta_{21} (I -\Theta_{11} )^{-1}\Delta_{11}(\tau) +w_2 \Delta_{22}(\tau)$.
Regarding part~\ref{convthirda} for $S_{22}$, for any $i \in S_{22}$ we have
\begin{equation*} 
\lim_{t \rightarrow \infty}  k_i(x(t)) =\lim_{t \rightarrow \infty}  \frac{\lambda^{t+1} (\alpha \nu_{2i} + \beta [(I -\Theta_{22} )^{-1} \Theta_{21} \nu_1]_i )}{\lambda^t (\alpha \nu_{2i} + \beta [(I -\Theta_{22} )^{-1} \Theta_{21} \nu_1]_i )} = \lambda. 
\end{equation*}
Regarding part~\ref{convthirdb} for $S_{22}$, notice that the theorem does not discuss the case with equal spectral radii.

Finally, statement~\ref{convforth} is proved utilizing previous statements. For any $i \in G_r(x(t))$ two cases exists. First, if $i$ belongs to a  closed-minded SCC, then $x_i(t) = x_i(\tau +1)$ for all $t > \tau$, and hence $k_i(x(t))$ does not exist. Second, if $i$ belongs to an open-minded SCC $S_k(x(t))$, then according to part~\ref{convthirda}, either $x_i(t) = \yf_i(t)$ or $k_i(x(t))$ converges to the spectral radius of the adjacency matrix of $S_k(x(t))$'s leader SCC. This spectral radius is proved in Section~\ref{blocks} to be strictly larger than zero and strictly smaller than one. In other words, there exists time $t_2$ such that for all $t \ge t_2$, $ 0 < k_i(x(t)) < 1$. Therefore, according to Remark~\ref{pseudoperstep}, $x(t)$ exhibits pseudo-stable behavior.
\end{proof}

\begin{figure}[h!]
  \centering
  \includegraphics[width=.47\textwidth,keepaspectratio]{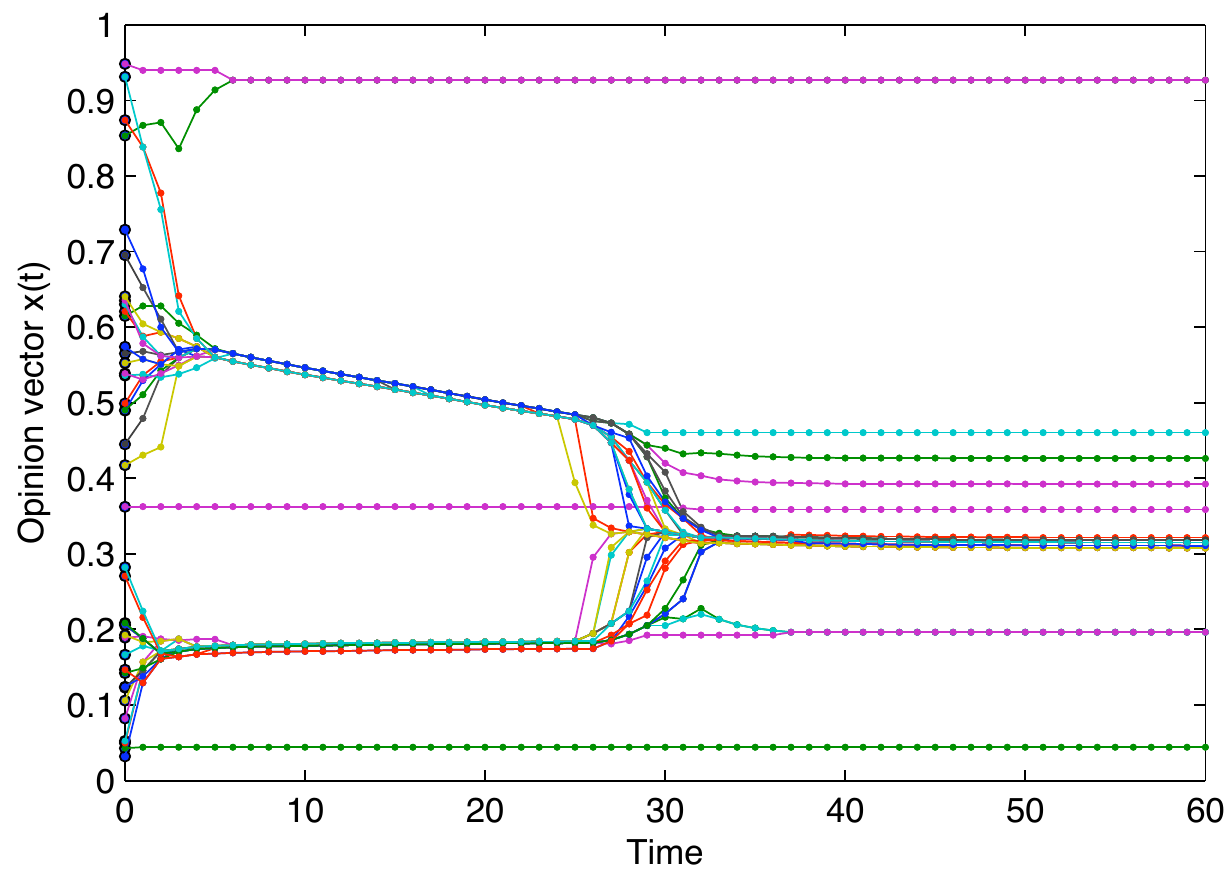}
  \includegraphics[width=.47\textwidth,keepaspectratio]{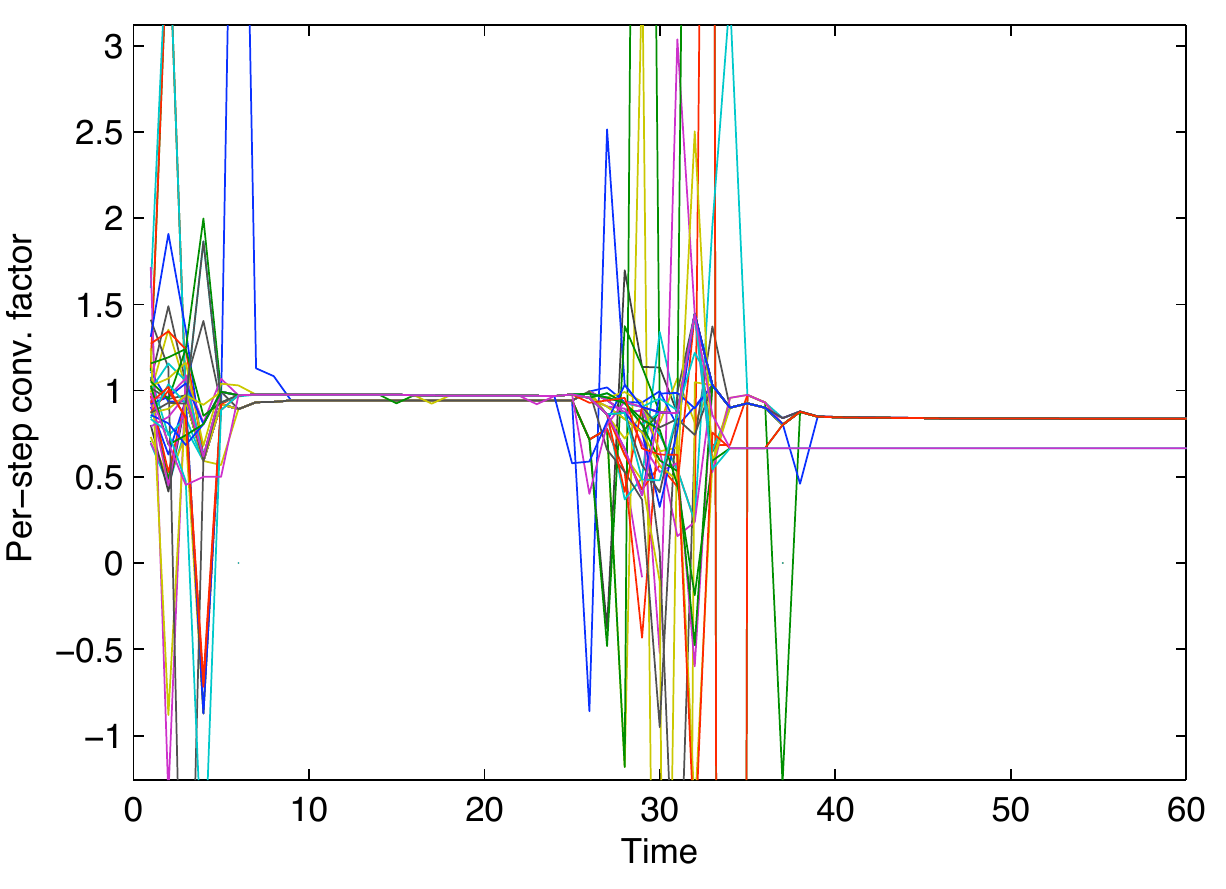}
  \caption{The trajectory of an SBC system (left) and the per-step
    convergence factor of its open-minded agents (right) are
    illustrated. The initial opinion vector and confidence bounds vector
    are generated randomly. This system satisfies the condition of
    Theorem~\ref{epsilonthm} at $t=50$. Moreover, since $x(50) \in
    \subscr{\mathcal{B}}{iet}(\yf(x(50)))$, the SBC digraph $G_r(x(t))$ is
    equal to $G_r(\yf(x(50)))$ for all $t\ge 50$. The digraph
    $G_r(\yf(x(50)))$ contains two open-minded SCC's, denoted by $S_1$ and
    $S_2$, while $S_2$ is a predecessor of $S_1$. The spectral radii of the
    adjacency matrices of $S_1$ and $S_2$ are equal to 0.6667 and 0.8381,
    respectively. Therefore, both $S_1$ and $S_2$ are their own leader
    SCC's, and by Theorem~\ref{convlemma} the per-step convergence factors
    of their agents converge to 0.6667 and 0.8381, respectively. The right
    plot verifies that the per-step convergence factors of all open-minded
    agents converge to those two values.}
\label{itsownleaderfig}
\end{figure}

\begin{figure}
  \centering
$\begin{array}{ccc}
  \includegraphics[width=2in,keepaspectratio]{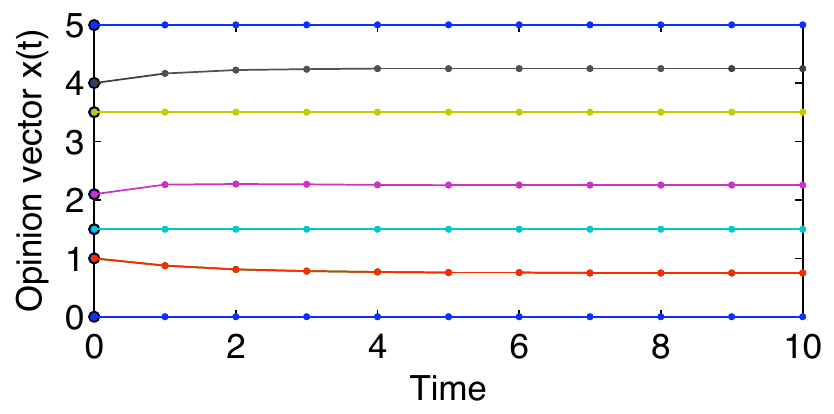}&
  \includegraphics[width=2in,keepaspectratio]{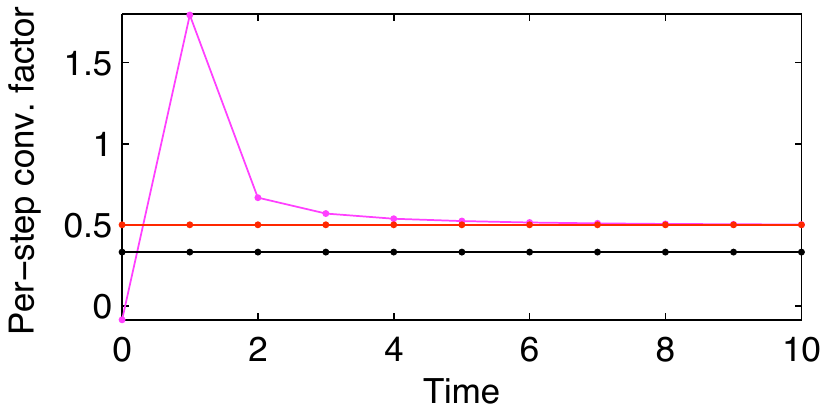} \\
\includegraphics[width=2in,keepaspectratio]{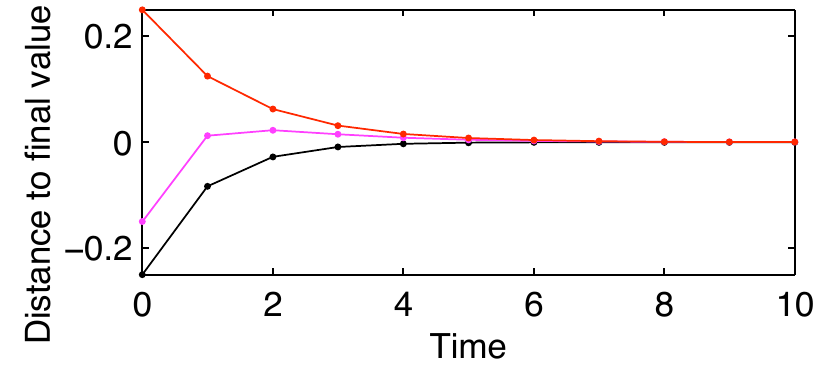}&
  \includegraphics[width=.8in,keepaspectratio]{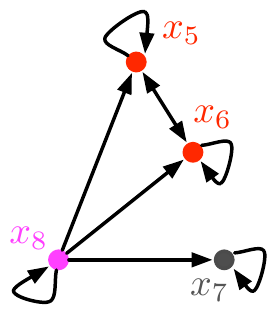} 
\end{array}$
\caption{An SBC trajectory $x(t)$ is
  plotted on the top left, the open-minded agents per-step convergence
  factors on the top right, the open-minded agents distances to their final
  values at constant topology $x_i(t) - x_i^*(x(t))$ on the bottom left,
  and the open-minded subgraph of $G_r(x(t))$ is illustrated on the bottom
  right.  This system is simulated with the initial vector $x(0) = [0
  \; 1.5 \; 3.5 \; 5 \; 1 \; 1 \; 4 \; 2.1]^T$ and confidence bounds $r =
  [0.01 \; 0.01 \; 0.01 \; 0.01 \; 1 \; 1 \; 1 \; 3]^T$. 
For all $t \ge 0$, the SBC digraph $G_r(x(t))$ remains unchanged
  and contains three open-minded SCC's: $\{ x_5, x_6 \}$, $\{ x_7 \}$, and
  $\{ x_8\}$. The spectral radii of the adjacency matrices of these SCC's
  are $0.5$, $0.333$, and $0.125$, respectively. The two SCC's $\{ x_5, x_6
  \}$ and $\{ x_7 \}$ are successors of $\{ x_8\}$, and based on their
  spectral radii, $\{ x_5, x_6\}$ is $\{ x_8\}$'s leader SCC.  We can see
  that the per-step convergence factor of $x_8$ converges to
  $0.5$. Furthermore, the sing of its direction of convergence toward the
  final value, i.e, the sign of $x_8(t) - x_8^*$, is the same as the
  leader's after $t=1$. These facts support Theorem~\ref{convlemma}. }
\label{differentdirfig}
\end{figure}

In Figures~\ref{itsownleaderfig} and \ref{differentdirfig}, we provide
numerical examples to facilitate the understanding of the conditions and
results of Theorem~\ref{convlemma}.

\section{Conclusion and Future Work}\label{conclusion}
This paper introduced a synchronized bounded influence (SBI) model of
opinion dynamics, which is similar to the heterogeneous bounded confidence
model introduced by Hegselmann and Krause, which we called synchronized
bounded confidence (SBC) model.  First, we conjectured that in both SBC and
SBI systems, for each trajectory there exists a finite time, after which
the topology of the interconnection network remains unchanged, hence, the
trajectory converges to a limiting opinion vector.  Second, we conjectured
that if a trajectory does not reach a fixed profile in finite time, then it
eventually shows a pseudo-stable behavior. We partly proved our first
conjecture, and the second conjecture is proved assuming that the first one
is true.  We designed a classification of agents that is employed in
computing the equilibria of the system.  We introduced the equi-topology
neighborhood and the invariant equi-topology neighborhood of the equilibria
of the system. Based on these neighborhoods, we derived sufficient
condition for both SBC and SBI systems to guarantee that the
interconnection topology remains unchanged for infinite time in a
trajectory, and therefore, the trajectory converges to a steady state.  In
our simulation results, it is observed that for uniformly randomly
generated initial opinion vector and bounds vector, the trajectories of
both systems eventually satisfy the mentioned sufficient condition with
probability one. However, the eventual convergence of every trajectory of
the SBC and SBI systems to a steady state is still an open problem.  Third,
we conjectured that, for uniformly randomly generated initial opinion
vector and bounds vector, the simulations of SBI systems converge in fewer
time steps and more often in finite time than SBC systems.  We derived
sufficient conditions for convergence in finite time for SBC and SBI
systems separately that intuitively explains our third conjecture.
Finally, we studied the trajectory of both SBC and SBI systems when they
update their opinions under fixed interconnection topology for infinite
time. We showed the existence of a leader group for each group of agents
that determines the follower's rate and direction of
convergence. 

The main future challenge is to prove that all SBC and SBI systems converge
to steady states. One approach is to prove that in each system, any
trajectory is eventually confined to the invariant equi-topology
neighborhood of an equilibrium opinion vector of the
system.
Moreover, the fact that the SBI systems converge in finite time more often
than the SBC systems might be explained by a probability analysis on the
topology of proximity digraphs.

\bibliographystyle{plain}
\bibliography{../../svn2/ref/alias,../../svn2/ref/FB,../../svn2/ref/Main}
% \bibliography{alias,FB,Main}

\end{document}